\documentclass[a4paper,twoside]{article}

\usepackage{epsf}
\usepackage{graphics}
\usepackage{latexsym}
\usepackage{amssymb}

\newcommand{\commentaar}[1]{{}}

\newcommand{\basis}[1]{\langle #1 \rangle}
\newcommand{\breukje}[2]{\mbox{\small $\frac{#1}{#2}$}}
\newcommand{\gen}[1]{\basis{#1}}
\newcommand{\halfje}{\breukje{1}{2}}
\newcommand{\diag}[1]{{\rm diag}(#1)}
\newcommand{\inprod}[2]{\langle #1, #2 \rangle}
\newcommand{\Ad}{{\rm Ad}}
\newcommand{\Gl}{{\rm \bf GL}}
\newcommand{\G}{{\rm \bf G}}
\newcommand{\Ker}{{\rm Ker}}
\newcommand{\Orb}{{\rm Orb}}
\newcommand{\Sl}{{\rm \bf Sl}}
\newcommand{\ad}{{\rm ad}}

\newcommand{\bfc}{{\rm \bf c}}
\newcommand{\bfi}{{\rm \bf i}}
\newcommand{\bfq}{{\rm \bf q}}
\newcommand{\bfs}{{\rm \bf s}}
\newcommand{\bft}{{\rm \bf t}}
\newcommand{\bfu}{{\rm \bf u}}
\newcommand{\defdas}{:=}
\newcommand{\eps}{\varepsilon}
\newcommand{\fC}{\mathbb{C}}
\newcommand{\fH}{\mathbb{H}}
\newcommand{\fR}{\mathbb{R}}
\newcommand{\fZ}{\mathbb{Z}}
\newcommand{\gl}{{\rm \bf gl}}
\newcommand{\g}{{\rm \bf g}}
\newcommand{\tB}{\tilde{B}}
\newcommand{\trace}{{\rm Trace}}
\renewcommand{\L}{{\cal L}}
\renewcommand{\sl}{{\rm \bf sl}}
\renewcommand{\sp}{{\rm \bf sp}}

\newcommand{\phm}{\hphantom{-}}

\newtheorem{theorem}{Theorem}[section]
\newtheorem{proposition}[theorem]{Proposition}
\newtheorem{lemma}[theorem]{Lemma}

\newtheorem{definition}{Definition}[section]
\newtheorem{remark}{Remark}[section]
\newtheorem{remarks}[remark]{Remarks}
\newtheorem{example}[remark]{Example}

\newcommand{\abc}{
\renewcommand{\theenumi}{\alph{enumi}}
\renewcommand{\labelenumi}{\theenumi)}
\itemsep 0pt
}

\makeatletter
\def\adots{\mathinner{\mkern1mu\raise\p@\vbox{\kern7\p@\hbox{.}}\mkern2mu
    \raise4\p@\hbox{.}\mkern2mu\raise7\p@\hbox{.}\mkern1mu}}
\def\ddots{\mathinner{\mkern1mu\raise8\p@\vbox{\kern7\p@\hbox{.}}\mkern2mu
    \raise4\p@\hbox{.}\mkern2mu\raise0\p@\hbox{.}\mkern1mu}}
\makeatother

\newcommand{\versie}{02/apr/05}

\title{Normal Forms and Unfoldings of Linear Systems in Eigenspaces of
(Anti)-Automorphisms of Order Two}
\author{I. Hoveijn\\
\small Langewoldlaan 2,\\
\small 9727 DD Groningen, Netherlands\\\\
J.S.W. Lamb\\
\small Department of Mathematics, Imperial College\\
\small London SW7 2BZ, UK\\\\
R.M. Roberts\\
\small Deparment of Mathematics \& Statistics, University of Surrey\\
\small Guildford GU2 7XH, UK}
\date{\versie}

\begin{document}
\maketitle
\begin{abstract}\noindent
In this article we classify normal forms and unfoldings of linear maps in
eigenspaces of (anti)-automorphisms of order two. Our main motivation is
provided by applications to linear systems of ordinary differential
equations, general and Hamiltonian, which have both time-preserving and
time-reversing symmetries. However the theory gives a uniform method to
obtain normal forms and unfoldings for a wide variety of linear
differential equations with additional structure. We give several examples
and include a discussion of the phenomenon of orbit splitting. As a
consequence of orbit splitting we observe \textit{passing} and
\textit{splitting} of eigenvalues in unfoldings.
\end{abstract}

\tableofcontents

%
\section{Introduction}\label{sec:intro}
Let $V$ be a finite dimensional real vector space and $\g$ an eigenspace of
an (anti)-automorphism $\gamma$ of order two of the Lie algebra
$\gl(V)$. Let $\G$ be a Lie subgroup of $\Gl(V)$ consisting of
\emph{structure preserving transformations} such that the action
\begin{displaymath}
L \mapsto gLg^{-1} \;\mbox{with $L \in \g$ and $g \in \G$},
\end{displaymath}
preserves $\g$. Then the $\G$-orbit of $L \in \g$ given by $\Orb_{\G}(L) =
\{gLg^{-1} \;|\; g \in \G\}$ is again a subset of $\g$. In this paper we
address the following two problems:
\begin{itemize}\parskip 0pt \itemsep 0pt
\item[i)] Classify all $\G$-orbits (normal forms) of elements $L$ in $\g$;
\item[ii)] Find the unfoldings of $L$ in $\g$.
\end{itemize}
We also briefly consider generalizations to abelian groups of
(anti)-automorphisms of order two.

%
\subsection*{Setting and Motivation}
Any local study of equilibrium points of vector fields starts with an
analysis of their linearizations. These are in one-to-one correspondence
with linear maps. This correspondence respects both the transformation
properties of linear vector fields under linear coordinate changes and
their Lie algebra structure.  Moreover dynamical systems theory is often
concerned with vector fields which preserve some structure. Well-known
examples are equivariant, reversible and Hamiltonian vector fields. The
linearizations of such vector fields preserve the same structure and the
spaces of structure preserving linear maps can be identified with
eigenspaces of (anti)-automorphisms of order two acting on the space of all
linear maps.

\begin{example}{\rm \label{ex:motiv1}
Consider reversible linear vector fields on $\fR^n$. Such a vector field is
determined by an infinitesimally reversible linear map $L$ satisfying
$RL=-LR$, where the linear \emph{structure map} $R$ satisfies $R^2=I$ and
$R \neq \pm I$.  We can also write this condition as $\phi_R(L) = -L$ where
the \emph{automorphism} $\phi_R$ is defined as $\phi_R(A) \defdas R^{-1}AR$
for all $A \in \gl(\fR^n)$. Thus $\g = \{A \in \gl(\fR^n) \;|\; \phi_R(A) = -A
\}$ is the $-1$ eigenspace of $\phi_R$. The structure preserving
transformation group $\G$ consists of $R$-equivariant maps $\G = \{g \in
\Gl(\fR^n)\;|\;gR=Rg\}$, that is elements of $\G$ map $\g$ into itself. See
section \ref{sec:eigauts} for a precise definition of the structure
preserving transformation group. \hfill $\rhd$}
\end{example}

\begin{example}{\rm \label{ex:motiv2}
Similarly a Hamiltonian linear vector field is determined by an
infinitesimally symplectic linear map $L$. Let $\omega$ be a symplectic
form, ie a non-degenerate skew symmetric bilinear form, on $\fR^{2n}$.  Let
$\inprod{\cdot}{\cdot}$ be the standard inner product on $\fR^{2n}$. Then
there is a \emph{structure map} $J$ satisfying $J^*=-J$ and $J^2=-I$ such
that $\omega(x,y) = \inprod{x}{Jy}$ for all $x,y \in \fR^{2n}$. An
infinitesimally symplectic map $L$ satisfies $\omega(Lx,y) = -\omega(x,Ly)$
or equivalently $\inprod{x}{L^*Jy} = -\inprod{x}{JL}$ for all $x,y \in
\fR^{2n}$. We write this condition as $\psi_J(L) = -L$, where $\psi_J$ is
the \emph{anti-automorphism} defined by $\psi_J(A) \defdas J^{-1}A^*J$, for
all $A \in \gl(\fR^{2n})$. Again $\g = \{A \in \gl(\fR^{2n}) \;|\;
\psi_J(A) = -A \}$ is the $-1$ eigenspace of $\psi_J$. In this case the
structure preserving transformation group $\G$ consists of maps that
preserve $\omega$, that is $\G = \{g \in \Gl(\fR^{2n})\;|\;\omega(gx,gy) =
\omega(x,y),\; {\rm for\; all}\; x,y \in \fR^{2n}\}$, which we can rephrase
using $J$ as $\G = \{g \in \Gl(\fR^{2n})\;|\; g^*Jg=J\}$. \hfill $\rhd$}
\end{example}

\begin{example}{\rm \label{ex:motiv3}
Combining the previous two examples, a reversible Hamiltonian linear vector
field is determined by an infinitesimally reversible symplectic linear map
$L$. Usually one requires that $R$ is an anti-symplectic map, and then
$\phi_R$ and $\psi_J$ commute. Thus infinitesimally reversible symplectic
maps on $\fR^{2n}$ are elements of the intersection of two eigenspaces
\begin{displaymath}
\{A\in\gl(\fR^{2n}) \;|\; \phi_R(A)=-A\} \cap
\{A\in\gl(\fR^{2n}) \;|\; \psi_J(A)=-A\},
\end{displaymath}
which, by virtue of the fact that $\phi_R$ and $\psi_J$ commute, is the
simultaneous eigenspace of $\phi_R$ and $\psi_J$. The structure preserving
transformation group for infinitesimally reversible symplectic maps is the
intersection of the transformation groups of Examples \ref{ex:motiv1} and
\ref{ex:motiv2}: $\G = \{g \in \Gl(\fR^{2n})\;|\; gR=Rg,\;g^*Jg=J\}$.
\hfill $\rhd$}
\end{example}

The main motivation for the theory developed in this paper is the normal
form and unfolding problem for linear reversible equivariant vector fields
in both the general and Hamiltonian cases. The spaces of such vector fields
have been described by Lamb \& Roberts \cite{lr} in the general case and
can be characterized as simultaneous eigenspaces of abelian groups of
(anti)-automorphisms of order two. The theory developed in this article
provides a uniform approach to all such problems. The authors plan to
report on applications of this theory to linear (Hamiltonian) reversible
equivariant vector fields in forthcoming publications.

In some cases, including Hamiltonian and equivariant vector fields, the
corresponding eigenspaces are Lie subalgebras of $\gl(V)$ and the normal
form and unfolding theory for maps in $\gl(V)$ (see Section
\ref{sec:nfglv}), carries over almost automatically.  However reversible
vector fields, for example, do not form a Lie subalgebra. This paper shows
that, despite this, analogous normal form and unfolding theories can been
developed.

Normal form and unfolding problems have a long history ranging from the
classical Jordan normal form to the more modern unfolding theory of Arnol'd
\cite{arn1}. We give a brief overview without trying to be
complete. Williamson \cite{wil} was the first to find normal forms for
infinitesimally symplectic maps. Later a more constructive approach was
given by Burgoyne \& Cushman \cite{bc1,bc2}. In this article we follow
their approach to a large extent. Unfoldings of infinitesimally
symplectic maps were independently given by Galin \cite{gal} and Ko\c{c}ak
\cite{koc}. For extensive studies of particular systems also see van der
Meer \cite{mee} and Cotter \cite{cot}. Normal forms and unfoldings of
infinitesimally reversible maps were first studied by Palmer \cite{pal} and
later by Sevryuk \cite{sev1} and Shih \cite{shi}. A particular example
where the linear part plays a crucial role can be found in Iooss
\cite{ioo}. Other contributions without a direct relation to dynamical
systems are Dempwolff \cite{dem}, Jacobson \cite{jac} for semi-linear maps
and Djukovic et al. \cite{dpwz} and Patera \& Rousseau \cite{pr2} for
subspaces of $\gl(V)$ which are not Lie algebras. Wiegman \cite{wie}
considers normal forms for maps over the quaternions. Studies of mixed
structures include Hoveijn \cite{hov} on infinitesimally reversible
symplectic maps and Melbourne \cite{mel} and Melbourne \& Dellnitz
\cite{md} on infinitesimally symplectic equivariant maps.

\begin{remark}{\rm
Note that the description using (anti)-automorphisms is not limited to
linear vector fields. In fact the latter are just the $1$-jets of
$C^{\infty}$-vector fields. The (anti)-automorphisms can equally well be
defined on $k$-jet spaces of vector fields, where they are still Lie
algebra (anti)-automorphisms of order two. The normalization procedures for
$C^{\infty}$-vector fields described in Broer et al. \cite{bdst} can be
combined with the ideas developed in this paper to give a corresponding
nonlinear normal form theory.\hfill $\rhd$}
\end{remark}

\begin{remark}{\rm
Another generalisation of the theory would be to consider a general
(compact) group $\Gamma$ of (anti)-automorphisms acting on $\gl(V)$
and classify normal forms and unfoldings of linear maps in
an isotypic component $\g$ of the action of $\Gamma$ on $\gl(V)$.}
\end{remark}

%
\subsection*{Main Results}
The main results of the paper are the Reduction Lemma \ref{lem:autreduc},
the Unfolding Lemma \ref{lem:autunfo} and the Orbit Splitting Theorem
\ref{the:split}. A formal statement of the Reduction Lemma requires some
technical notation, but it may be informally summarised as:

\textbf{Reduction Lemma}\\
\emph{The normal form of a linear map $L$ in an eigenspace of an
(anti)-automorphism is determined by the semi-simple part of $L$ on a
reduced space.}

This lemma greatly simplifies the problem of finding normal forms, because
the actual computations are limited to low dimensional spaces. It is
essential for the Reduction Lemma that the Jordan-Chevalley decomposition
holds in the eigenspaces of an (anti)-automorphism. Using the notation
introduced above, the Unfolding Lemma reads as follows.

\textbf{Unfolding Lemma}\\
\emph{Let $\gamma$ be an (anti)-automorphism of finite order and let $L \in
\g$. Then the restriction of the $\Gl(V)$-centralizer unfolding of $L \in
\gl(V)$ to $\g$ is equivalent to the $\G$-centralizer unfolding in $\g$.}

This means that we do not need to find a new way of computing unfoldings in
a subset of $\gl(V)$ with a smaller structure preserving transformation
group $\G$. We simply use a version of the existing Arnol'd or centralizer
unfolding, see Lemma \ref{lem:unfo}. As an alternative one might use the
representation theory of $\sl(2)$ to find unfoldings, see Ko\c{c}ak
\cite{koc} or Cushman \& Sanders \cite{cs}.

Orbit splitting is a well known phenomenon for infinitesimally symplectic
maps. If such a map has a pair of double imaginary eigenvalues then there
are two inequivalent normal forms. They may be distinguished by
\emph{signs}, see Example \ref{ex:sympl}. In general the $\Gl(V)$-orbit of
a map $L \in \g$ may intersect $\g$ in several $\G$-orbits of $L$. The
Orbit Splitting Theorem states that there are at most two such orbits.

\textbf{Orbit Splitting Theorem}\\
\emph{The intersection of the $\Gl(V)$-orbit of $L$ in $\gl(V)$ and $\g$
consists of at most two $\G$-orbits.}

In general inequivalent $\G$-orbits have different unfoldings, which may
give rise to \emph{passing} and \emph{splitting} of eigenvalues when
parameters are varied, see Section \ref{sec:unfos} for details.

%
\subsection*{Organization}
The remainder of this article is organized as follows. In Section
\ref{sec:nfglv} we review the theory for normal forms and unfoldings in
$\gl(V)$. We use this as a starting point for finding normal forms and
unfoldings in the eigenspace of an (anti)-automorphism in Section
\ref{sec:nfeig}. In Section \ref{sec:nfunfo} we apply the results of Section
\ref{sec:nfeig} to present normal forms and unfoldings in eigenspaces of
(anti)-automorphisms of order two.
Finally in Section \ref{sec:gens} we generalize our results to abelian
groups of (anti)-automorphisms of order two. We also suggest some further
possible generalizations.

%
\section{Normal Forms and Unfoldings in $\gl(V)$}\label{sec:nfglv}
A linear differential equation is given by $\dot{x} = Ax$ where $A \in
\gl(V)$. A coordinate change $y = gx$, with $g \in \Gl(V)$, transforms
this to $\dot{y} = gAg^{-1}y$. Thus
linear vector fields transform as linear maps. We therefore identify the
space of linear vector fields on $V$ with $\gl(V)$. Here we review the
normal form and unfolding theory for linear maps in $\gl(V)$ in Sections
\ref{sec:nufglv} and \ref{sec:unfglv} respectively.

%
\subsection{Normal Forms}\label{sec:nufglv}
Let $V$ be a finite dimensional real vector space. Then $\gl(V)$ is the Lie
algebra of all linear maps from $V$ to itself. The Lie group $\Gl(V)$ is
the group of all invertible linear transformations from $V$ to itself. The
action of $\Gl(V)$ on $\gl(V)$ is given by the \emph{adjoint action}, that
is, by \emph{similarity transformations}:
\begin{displaymath}
\Ad_g : L \mapsto gLg^{-1}
\end{displaymath}
The $\Gl(V)$-orbits
\begin{displaymath}
\Orb_{\Gl(V)}(L) \ =\ \{gLg^{-1} \;|\; g \in \Gl(V)\}
\end{displaymath}
of the adjoint action are precisely the equivalence classes we
are interested in classifying. From now on we will use the word `orbits'
only. It is well known that for the adjoint action of $\Gl(V)$ the orbit of
$L$ in $\gl(V)$ is determined by two invariants: the eigenvalues and Jordan
structure of $L$. The \emph{Jordan-Chevalley decomposition},
\emph{Reduction Lemma} and \emph{Reconstruction Lemma} described below
formalize this fact.

%
\subsubsection*{Jordan-Chevalley Decomposition}\label{sec:jcdglv}

The Jordan-Chevalley decomposition splits a linear map $L$ into the sum of
its semi-simple and nilpotent parts. In order to define semi-simple we need
to work over the complex numbers, so in this section we assume that $L$ is
defined on a complexified space $V$. In Theorem \ref{the:vlambda} we
translate our results for a real space $V$. A map $S$ is called
\emph{semi-simple} if the algebraic and geometric multiplicity of each of
its eigenvalues are equal. A map $N$ is called \emph{nilpotent} if $N^n=0$
for some integer $n$. The least such integer is called the \emph{height} of
$N$.
\begin{theorem}
\label{the:jcd}
\textbf{\emph{(Jordan-Chevalley decomposition)}}\\ For each $L \in \gl(V)$
there exist a unique semi-simple $S\in \gl(V)$ and a unique nilpotent $N
\in \gl(V)$ such that $[S,N]=0$ and $L = S + N$.
\end{theorem}
The eigenvalues of $L$ are determined by the semi-simple part $S$ while the
nilpotent part $N$ determines its Jordan structure.
The Jordan-Chevalley decomposition of a linear map is $\Ad_g$-equivariant
and so is a property of the $\Gl(V)$-orbit rather than the individual map,
see Humphreys \cite{hum}.
The Jordan-Chevalley decomposition almost automatically extends to all the
classical Lie algebras, see \cite{hum,var}. We shall see in
Section \ref{sec:nfeig} that it also extends to eigenspaces of Lie algebra
(anti)-automorphisms acting on $\gl(V)$.

%
\subsubsection*{Reduction}\label{sec:reducglv}
The Reduction Lemma, due to Burgoyne \& Cushman \cite{bc2,bc3}, exploits
the Jordan-Chevalley decomposition to simplify the normal form
problem for linear maps. It formalizes the classical Jordan normal form
algorithm and works in all classical Lie algebras. In the first part of
this section we will work on the complexified space $V$.

Let $L \in \gl(V)$ be a linear map with Jordan-Chevalley decomposition
$L=S+N$. An $L$-invariant subspace of $V$ is said to be
\emph{indecomposable} if it has no proper $L$-invariant subspaces. The
restriction of $L$ to an indecomposable $L$-invariant subspace has a unique
eigenvalue and such a subspace is a generalized eigenspace for that
eigenvalue. The space $V$ decomposes as a direct sum of indecomposable
$L$-invariant subspaces $V = \oplus_{\lambda} V_{\lambda}$ which is unique
up to permutations unless there are two or more $V_{\lambda}$ with equal
eigenvalues \emph{and} equal heights. Moreover the restrictions of $S$ and
$N$ to the indecomposable subspaces are the semi-simple and nilpotent parts
of the restrictions of $L$. The characteristic polynomial of $L$ factors
over the indecomposable $L$-invariant subspaces. The height $n$ of the
restriction of $N$ to such a space can be computed from the characteristic
polynomial of the restriction of $L$, see Burgoyne \& Cushman \cite{bc1}.
\begin{lemma}
\label{lem:reduc}
\textbf{\emph{(Reduction Lemma)}}\\ Let $V_{\lambda}$ be an indecomposable
$L$-invariant space with eigenvalue $\lambda$. Assume that the restriction
of $N$ to $V_\lambda$ has height $n$. Then there is an $S$-invariant
complement $W_{\lambda}$ of $NV_{\lambda}$ in $V_{\lambda}$ such that
$V_{\lambda}=W_{\lambda} \oplus NV_{\lambda} = W_{\lambda} \oplus
NW_{\lambda} \oplus \cdots \oplus N^{n-1} W_{\lambda}$. For each $j =
0\ldots n-1$ we have $\dim N^jW_{\lambda}=\dim W_{\lambda}$.  The
restriction of $L$ to $V_{\lambda}$ is determined up to similarity by the
restriction of $S$ to $W_{\lambda}$.
\end{lemma}

Thus if we wish to classify a linear map $L$ we only have to classify its
semi-simple part $S$ by Lemma \ref{lem:reduc}. From now on we only work
over the real numbers so $V$ is again a real space. With a slight abuse of
notation we write again $V_{\lambda}$ even if $\lambda$ is complex. Theorem
\ref{the:vlambda} relates the real spaces to the complex ones. On many
occasions we distinguish four cases. If $\lambda$ is zero or
$\lambda=\alpha$ is real we write $V_0$ and $V_{\alpha}$ for the
generalized eigenspaces. If $\lambda = \pm i\beta$ is purely imaginary or
$\lambda = \alpha \pm i\beta$ is complex we write $V_{\pm i\beta}$ and
$V_{\alpha \pm i\beta}$ for the real generalized eigenspaces.

\begin{theorem}\label{the:vlambda}
Let $L=S+N$ be the Jordan-Chevalley decomposition of $L$ and let
$V_{\lambda}$ be an indecomposable $L$-invariant subspace with eigenvalue
$\lambda$. Assume that the restriction of $N$ to $V_{\lambda}$ has height
$n$.
\begin{itemize}\itemsep 0pt
\item[a)] {\bf Real eigenvalues.}
If $\lambda=\alpha \in \fR$ then $(L-\alpha)^n = 0$. In this case $\dim
W_{\alpha}=1$, $\dim V_{\alpha}=n$ and for all $e \in V_{\alpha}$ we have
$Se=\alpha e$.
\item[b)] {\bf Complex eigenvalues.}
If $\lambda=\alpha \pm i\beta$ with $\alpha,\beta \in \fR$ and
$\beta \neq 0$ then $((L-\alpha)^2+\beta^2)^n = 0$. In this case $\dim
W_{\alpha \pm i\beta}=2$, $\dim(V_{\alpha \pm i\beta})=2n$ and for all $v
\in V_{\alpha \pm i\beta}$ we have $(S-\alpha)^2v=-\beta^2 v$. For any $e
\in W_{\alpha \pm i\beta}$ let $f = \frac{1}{\beta}(S-\alpha)e$, then
$\basis{e,f}$ is a basis of $W_{\alpha \pm i\beta}$.
\end{itemize}
\end{theorem}
Thus the restriction of the semi-simple part $S$ to a subspace
$W_{\lambda}$ always has normal form
$$
\left( \begin{array}{c}
               \alpha
        \end{array} \right)
\quad \quad \mbox{or} \quad \quad
\left( \begin{array}{cc}
               \alpha  &-\beta \\
               \beta   & \alpha
        \end{array} \right),
$$
depending on whether $\lambda$ is real or complex, respectively.
Note that in case of complex eigenvalues we can always find a basis such
that $\beta > 0$.

%
\subsubsection*{Reconstruction}\label{sec:reconglv}
Suppose we are given a linear map $L\in \gl(V)$ with Jordan-Chevalley
decomposition $L=S+N$ and an indecomposable $L$-invariant space
$V_{\lambda}$. From Lemma \ref{lem:reduc} we know that there is an
$S$-invariant complement $W_{\lambda}$ of $NV_{\lambda}$ in $V_{\lambda}$
such that $V_{\lambda}= W_{\lambda} \oplus NW_{\lambda} \oplus \cdots
\oplus N^{n-1}W_{\lambda}$, where $n$ is the height of $N$ on $V_{\lambda}$
and $\dim N^jW_\lambda = \dim W_\lambda$. In Theorem \ref{the:vlambda} we
gave normal forms for the restriction of $S$ to $W_{\lambda}$. To find
normal forms for $L$ on $V_{\lambda}$ we start with the basis of
$W_{\lambda}$ used for the normal form of $S$ on $W_{\lambda}$ in Theorem
\ref{the:vlambda} and apply $N$ to this $n-1$ times to generate a basis for
$V_{\lambda}$.
\begin{lemma}
\label{lem:recons}
\textbf{\emph{(Reconstruction Lemma)}}\\ Let $L=S+N$ be the
Jordan-Chevalley decomposition of $L$ and $V_{\lambda}$ be an
indecomposable $L$-invariant subspace.  Let $W_\lambda$ be an $S$-invariant
complement to $NV_\lambda$ in $V_\lambda$ and assume that $N$ has height
$n$ on $V_\lambda$.  Then, if $\basis{e_1\ldots,e_m}$ is a basis for
$W_\lambda$, the set $\basis{e_1\ldots,e_m, Ne_1\ldots, Ne_m,
N^{n-1}e_1\ldots,N^{n-1}e_m}$ is a basis for $V_\lambda$.
\end{lemma}

The corresponding normal forms for indecomposable $L$-invariant subspaces
with real and complex eigenvalues are the familiar real Jordan blocks
\begin{displaymath}
\left(\begin{array}{cccc}
\alpha &&&\\
1&\ddots&&\\
&\ddots&\ddots&\\
&&1&\alpha
\end{array}\right)
\quad\quad \mbox{and} \quad\quad
\left(\begin{array}{rrrrrrrr}
\alpha &-\beta &&&&&&\\
\beta &\alpha  &&&&&&\\
1& 0& \alpha &-\beta &&&&\\
0& 1& \beta &\alpha  &&&&\\
&&\ddots&&\ddots&&&\\
&&&\ddots&&\ddots&&\\
&&&&1&0&\alpha&-\beta\\
&&&&0&1&\beta&\alpha
\end{array}\right)
\end{displaymath}
with respect to the bases $\basis{e,Ne,\ldots,N^{n-1}e}$ and
$\basis{e,f,Ne,Nf,\ldots,N^{n-1}e,N^{n-1}f}$, respectively. Here $e$ and
$f$ are vectors as in Theorem \ref{the:vlambda}.

%
\subsection{Unfoldings}\label{sec:unfglv}
A general theory for \emph{unfoldings} or \emph{deformations} of maps in
$\gl(V)$ is given by Arnol'd \cite{arn1}. Here we use the Reduction
Lemma to describe this theory, see also \cite{hov}. Note that
$\Gl(V)$-orbits are smooth submanifolds of $\gl(V)$, see Bredon \cite{bre}.
\begin{definition}
A smooth map $\L: \fR^p \rightarrow \gl(V): \mu \mapsto \L(\mu)$ with
$\L(0)=L$ is called an \emph{unfolding} or a \emph{deformation} of $L$. If
$\L$ is transverse to the $\Gl(V)$-orbit through $L$ at $L$, then it is
said to be \emph{versal}.
\end{definition}
From now on we will only use the word `unfolding'. We are especially
interested in unfoldings of a map $L$ having a minimum number of parameters
but still parametrising a section transverse to the $\Gl(V)$-orbit through
$L_{0}$. We therefore make the following definition, see Arnol'd
\cite{arn1}.
\begin{definition}
Two unfoldings $A(\mu)$ and $B(\mu)$ of $L$ are called \emph{equivalent} if
they are similar as families of linear maps. This means that there is a
smooth family of transformation $g(\mu) \in \Gl(V)$ such that
$g(\mu)A(\mu)g(\mu)^{-1} = B(\mu)$ for all $\mu \in \fR^p$.  An unfolding
$\L$ of $L$ is called \emph{miniversal} if (a) for every other unfolding $A:
\fR^q \rightarrow \gl(V)$ of $L$ there exists a smooth map $\chi: \fR^q
\rightarrow \fR^p$ such that $A$ is equivalent to $\L \circ \chi$, and (b)
$\L$ has the minimal number of parameters possible for unfoldings with this
property.
\end{definition}
The number of parameters for a miniversal unfolding is equal to the
codimension of the $\Gl(V)$-orbit through $L$ and so is called the
\emph{codimension} of $L$. Geometrically, the image of a miniversal
unfolding of $L$ is a submanifold of $\gl(V)$ whose tangent space at
$L$ is a complement to the tangent space at $L$ of the $\Gl(V)$-orbit
through $L$.

Arnol'd \cite{arn1} showed that miniversal unfoldings can be obtained by
taking orthogonal complements to tangent spaces of $\Gl(V)$-orbits with
respect to the inner product $\inprod{A}{B} = \trace(A^*B)$ on $\gl(V)$.
Some computations show that the \emph{centralizer} or \emph{Arnol'd
unfolding} is:
\begin{displaymath}
\{L+M^* \;|\; M \in \Ker(\ad_L) \},
\end{displaymath}
where $\ad_L(M) = [M,L] = ML-LM$. By applying the adjoint action of
$\Gl(V)$ on $\gl(V)$ to this unfolding we can obtain an unfolding at any
other point on the $\Gl(V)$-orbit through $L$. Transversality and
miniversality are preserved by this transformation, but orthogonality will
usually be lost. Another way to find unfoldings is the use of
representation theory of $\sl(2)$, see Ko\c{c}ak \cite{koc} or Cushman \&
Sanders \cite{cs}.

We now give a more detailed description of the centralizer unfolding using
the Reduction Lemma. The idea is to start with the centralizer unfolding of
the semi-simple part of $L$ and then reconstruct the unfolding of $L$ in
several steps. An advantage of this approach is that it depends only on
the Jordan-Chevalley decomposition and not on a particular normal form for
$L$. Furthermore we only have to compute matrices which commute with a
semi-simple matrix on a low dimensional space.

Assume that $L=S+N$ has only one eigenvalue on $V$ and let $V=V_1 \oplus
\cdots \oplus V_p$, where each $V_i$ is an indecomposable $L$-invariant
subspace.  For each $V_i$ we have the decomposition $V_i=W_i \oplus NW_i
\oplus \cdots \oplus N^{n_i-1}W_i$, where $W_i$ is an $S$-invariant
complement to $NV_i$ in $V_i$. Assume that the heights $n_i$ of $N$
restricted to $V_i$ satisfy $n_1 \geq \cdots \geq n_p$.  Let $W = W_1
\oplus \cdots \oplus W_p$ and let $\basis{e_{i,1},\ldots, e_{i,q}}$ be a
basis for $W_i$.  Note that the $W_i$ all have the same dimension $\dim W_i
= q$, so we may assume that their bases have been chosen such that the
matrices of $S_i = S|_{W_i}$ are equal.

To find the centralizer unfolding of $L$ first choose a basis for
$\bfu_{W_1} = \{M : W_1 \rightarrow W_1 \;|\; \ad_S(M) = 0 \}$.  Only in
this step is it necessary to compute commuting matrices. The next step is to
extend the basis for $\bfu_{W_1}$ to a basis for $\bfu_W = \{M : W
\rightarrow W \;|\; \ad_S(M) = 0 \}$. Then we extend this basis to the set
of maps on $V$ which commute with $N$. The final step is to construct a
basis for $\bfu_V = \{M : V \rightarrow V \;|\; \ad_L(M) = 0 \}$. We make
this more precise in the following Unfolding Lemma, which is a corrected
version of that in \cite{hov}.
\begin{lemma}
\label{lem:unfo}
\textbf{\emph{(Unfolding Lemma)}}
\begin{enumerate}\itemsep 0pt
\item Construct a basis $\basis{B^{(1)},\ldots,B^{(r)}}$ of
$\bfu_{W_1}$.
\item Extend the basis of $\bfu_{W_1}$ to a basis of $\bfu_W$ by
defining $B^{(k)}_{ij}:W \rightarrow W$ by:
\begin{displaymath}
B^{(k)}_{ij}\ =\ \left\{\begin{array}{ll} B^{(k)}& :\ W_j \rightarrow W_i\\
0 & :\ W_{j'} \rightarrow W_{i'}, \; (i',j')\neq(i,j)
\end{array}\right.
\end{displaymath}
Then $\{B^{(k)}_{ij}\;|\;i,j=1,\ldots,p \;{\rm and}\;k=1,\ldots,r\}$
is a basis of $\bfu_W$.
\item Extend the basis of $\bfu_W$ to a basis of $\bfu_V$ by defining
$\tB^{(k)}_{ij}:V\rightarrow V$ by:
\begin{displaymath}
\tB^{(k)}_{ij} N^m e_{j,l}\ =\
\left\{\begin{array}{ll}
N^m B^{(k)}_{ij} e_{j,l} & j \geq i\\
N^{m+n_i-n_j} B^{(k)}_{ij} e_{j,l} & j < i
\end{array}\right.
\end{displaymath}
for $l=1,\ldots,q$ and $m=0,\ldots,n_j-1$. Then the $N^n \tB^{(k)}_{ij}$
for $0 \leq n < \min(n_i,n_j)$ form a basis of $\bfu_V$.
\item Let $\basis{B_1,\ldots,B_d}$ be a basis of $\bfu_V$, then $L(\mu) = L
+ \sum_{i=1}^d \mu_{i} B_i^*$ is a miniversal deformation of $L$.
\end{enumerate}
\end{lemma}
\commentaar{bundle codim versus orbit codim ?  tja...}
\textbf{Proof.} The construction in the proposition is a formalization of
the construction by Arnol'd \cite{arn1}, which in turn is based on
Gantmacher's construction in \cite{gan}. We restrict ourselves to the case
that $L$ has a real eigenvalue. First we prove that the $B_1,\ldots,B_d$
are linearly independent, then we prove that they span $\bfu_V$. By
construction the $B^{(1)},\ldots,B^{(r)}$ are linearly independent. Then
the $B^{(k)}_{ij}$ are also linearly independent since they map $W_j$ to
$W_i$. In step 3) the $B^{(k)}_{ij}$ are only extended to $V$ so the
$\tB^{(k)}_{ij}$ are still linearly independent. (Here we have to take into
account that the height of $N$ on $W_j$ may be smaller than the height of
$N$ on $W_i$ so that $N^{m+n_i-n_j}=0$ as soon as $N^m=0$.) The $N^n
\tB^{(k)}_{ij}$ are linearly independent because they map $W_j$ to
$N^nW_i$. The number of $N^n \tB^{(k)}_{ij}$ is equal to $\sum_{i,j}
\min(n_i,n_j) = \sum_i (2i-1)n_i$ which is equal to the number of
parameters in the Arnol'd unfolding, so the $N^n \tB^{(k)}_{ij}$ span
$\bfu_V$. Thus the $N^n \tB^{(k)}_{ij}$ form a basis of $\bfu_V$.\hfill
$\Box$


\section{Normal Forms and Unfoldings in an Eigenspace of an
(Anti)-Automorphism}\label{sec:nfeig}
We outline a general theory for maps in eigenspaces of Lie algebra
(anti)-automorphisms of order two acting on $\gl(V)$ and show that a
$\Gl(V)$-orbit in $\gl(V)$ can split into at most two $\G$-orbits when
intersected with such an eigenspace. Here $\G$ is the structure preserving
transformation group.

To characterize normal forms in Section \ref{sec:nufglv} and find their
unfoldings in Section \ref{sec:unfglv} we made essential use of the
following facts. First, the equivalence classes are the orbits of the
adjoint action of a transformation group. Second, the Jordan-Chevalley
decomposition leads to a reduction of the normal form and unfolding
problems to semi-simple maps. Third, the Lie algebra of the transformation
group enables us to characterize tangent spaces to orbits and hence to find
miniversal unfoldings as complements.

In this section we show that the eigenspaces of (anti)-automorphisms have
all the Lie algebraic structure that is needed to generalize these facts.
Appropriate \emph{structure preserving} transformation Lie groups can be
defined and the Jordan-Chevalley decomposition still holds. The latter
leads to a Reduction Lemma which can be used to simplify the computation of
normal forms and miniversal unfoldings.

%
\subsection{(Anti)-Automorphisms}\label{sec:aaut}
We begin by describing some properties of (anti)-automorphisms of
$\gl(V)$. Where it is convenient we identify $V$ with $\fR^n$ and hence
$\gl(V)$ with $\gl(n,\fR)$.
\begin{definition}
\label{def:aut}
A linear map $\gamma: \gl(V) \rightarrow \gl(V)$ is an {\rm automorphism}
of $\gl(V)$ if $\gamma(AB) = \gamma(A) \gamma(B)$ and an {\rm
anti-automorphism} if $\gamma(AB) = \gamma(B) \gamma(A)$.
\end{definition}
We will generally denote an automorphism by $\phi$, an anti-automorphism by
$\psi$ and will write $\gamma$ if we do not want to make a distinction. Let
\begin{displaymath}
\sigma(\gamma) \ =\ \left\{\begin{array}{l}
\phm 1, \mbox{ if $\gamma$ is an automorphism}\\
-1, \mbox{ if $\gamma$ is an anti-automorphism.}
\end{array}\right.
\end{displaymath}
Note that for any $\gamma$ the Lie bracket $[A,B] = AB-BA$ on $\gl(V)$
satisfies
\begin{displaymath}
\gamma\left(\left[A,B\right]\right) \ =\
\sigma(\gamma) \left[ \gamma(A), \gamma(B) \right].
\end{displaymath}
Thus $\gamma$ is a Lie algebra automorphism or anti-automorphism.

The next proposition describes all (anti)-automorphisms of $\gl(V) \cong
\gl(n,\fR)$ and shows that they have associated \emph{structure maps} $\bfs
\in \Gl(n,\fR)$. A proof based on the symmetries of the Dynkin diagram can
be found in Freudenthal \& de Vries \cite{fv}.
\begin{proposition}\label{pro:inout}
Every automorphism of $\gl(n,\fR)$ has the form $\phi_{\bfs}(L) = \bfs^{-1}
L \bfs$, $L \in \gl(n,\fR)$, for some $\bfs \in \Gl(n,\fR)$.  The
anti-automorphisms of $\gl(n,\fR)$ are all of the form $\psi_{\bfs} =
\phi_{\bfs} \circ \psi_I$ where $\psi_I(L) = L^*$.
\end{proposition}
Here $L^*$ is defined as $\inprod{x}{L^*y} \defdas \inprod{Lx}{y}$, for all
$x,y \in \fR^n$ and $\inprod{\cdot}{\cdot}$ is an inner product on
$\fR^n$. In the next lemma we collect some simple properties of an
(anti)-automorphism of order two.
\begin{lemma}\label{lem:gprops}
Let $\gamma_{\bfs}$ be an (anti)-automorphism of order two on $\gl(V)$ with
structure map $\bfs \in \Gl(V)$. Then:
\begin{enumerate}\abc\parskip 0pt \itemsep 0pt
\item $\gamma_{\alpha \bfs} =\gamma_{\bfs}$ for all $\alpha \in \fR$;
\item The eigenvalues of $\gamma_{\bfs}$ are $\pm 1$;
\item $\gamma_{\bfs}$ is semi-simple.
\end{enumerate}
\end{lemma}

Note that a) holds for any (anti)-automorphism and c) is true for
(anti)-automorphisms of finite order. The next proposition lists some
properties of structure maps.
\begin{proposition}\label{pro:sprops}
Let $\bfs$ be a structure map associated to an (anti)-automorphism
$\gamma_{\bfs}$ of order two. Then $\bfs$ has the following properties:
\begin{enumerate}\abc\parskip 0pt \itemsep 0pt
\item $\det \bfs = \pm 1$;
\item $\bfs$ is orthogonal;
\item $\bfs^2 = \pm I$;
\item There exists a basis such that the matrix of $\bfs$ is given by $\bfs =
\left(\begin{array}{rr}I_p&0\\0&-I_q\end{array}\right)$ or $\bfs =
\left(\begin{array}{rr}0&-I_n\\I_n&0\end{array}\right).$
\end{enumerate}
\end{proposition}
Here a) is true in general because it is a consequence of property a) in
Lemma \ref{lem:gprops} and b) holds for any (anti)-automorphism of finite
order.

\textbf{Proof of Proposition \ref{pro:sprops}.} Using $\gamma_{\alpha \bfs}
=\gamma_{\bfs}$ for all $\alpha \in \fR$, we can scale $\bfs$ such that
$\det \bfs = \pm 1$. For the remaining parts we distinguish between
automorphisms and anti-automorphisms. We use the fact that $\bfs$
transforms as $\bfs \mapsto g \bfs g^{-1}$ or $\bfs \mapsto g \bfs g^*$
under a coordinate transformation $g$ depending on whether $\bfs$ is
associated to an automorphism or an anti-automorphism, see Section
\ref{sec:eigauts}. Let $\gamma_{\bfs}=\phi_{\bfs}$ be an automorphism. Then
$\phi_{\bfs}^2=I$ immediately implies $\bfs^2 = \pm I$. This in turn
implies that $\bfs$ has eigenvalues $\pm 1$ or $\pm i$ and that $\bfs$ is
semi-simple, which yields the normal forms in d). Thus after a suitable
transformation $\bfs$ is orthogonal. Let $\gamma_{\bfs}=\psi_{\bfs}$ be an
anti-automorphism. Then $\psi_{\bfs}^2=I$ implies $\bfs^*=\pm \bfs$ and so
is semi-simple. Thus $\bfs$ has either real or purely imaginary
eigenvalues. Restrict $\bfs$ to an indecomposable $\bfs$-invariant space,
then by applying the map $g=\rho I$ the eigenvalues of $\bfs$ are scaled to
$\pm 1$ or $\pm i$. This transformed $\bfs$ is orthogonal and moreover
satisfies $\bfs^2=\pm I$. Thus it takes one of the forms in d) with respect
to a suitable basis. \hfill $\Box$

%
\subsection{Eigenspaces of (Anti)-Automorphisms}\label{sec:eigauts}
The eigenspace of an (anti)-automorphism $\gamma$ corresponding to an
eigenvalue $\mu$ is denoted by
\begin{displaymath}
\gl_{\mu}(V) \ =\
\{L \in \gl(V) \;|\; \gamma(L) = \mu L\}.
\end{displaymath}
Here we only consider (anti)-automorphisms of order two so that $\gl_{\mu}(V)$
is again a real space. The next lemma gives some basic properties of the
eigenspaces of $\gamma$.
\begin{lemma}\label{lem:almod}
Let $\gamma$ be an (anti)-automorphism of order two on $\gl(V)$. Then:
\begin{enumerate}\abc\itemsep 0pt
\item $\gl_\mu(V)$ is a Lie subalgebra of $\gl_\mu(V)$ if and only if $\mu
= \sigma(\gamma)$;
\item $\gl_\mu(V)$ is a Lie submodule of $\gl(V)$ over the Lie subalgebra
$\gl_{\sigma(\gamma)}(V)$;
\item $\gl(V)$ splits as a sum of eigenspaces, $\gl(V) = \gl_1(V) \oplus
\gl_{-1}(V)$.
\end{enumerate}
\end{lemma}
Note that a) and b) hold for any (anti)-automorphisms of finite order.

\textbf{Proof of Lemma \ref{lem:almod}.} Let $A \in \gl_{\mu}(V)$ and $B
\in \gl_{\nu}(V)$. Then $\gamma([A,B]) = \sigma(\gamma)[\gamma(A),
\gamma(B)] = \sigma(\gamma) \mu \nu [A,B]$. Since $\gamma$ is of order two
it is semi-simple and its eigenvalues are $\pm 1$, and so $\gl(V)$
splits as in c).\hfill $\Box$

%
\subsubsection*{Jordan-Chevalley Decomposition}\label{sec:autjcd}
The Jordan-Chevalley Decomposition \ref{the:jcd} holds in any Lie
subalgebra $\g$ of $\gl(V)$, see \cite{hum}. The next proposition states
that it holds in any eigenspace of an (anti)-automorphism of $\gl(V)$. Note
that $\gl_{\mu}(V)$ need not be real for the proposition to hold.
\begin{proposition}\label{pro:autjcd}
If $L \in \gl_{\mu}(V)$ and $L=S+N$ with $S$ semi-simple, $N$
nilpotent and $[S,N]=0$, then both $S$ and $N$ are elements of
$\gl_{\mu}(V)$.
\end{proposition}

{\bf Proof.} Let $L \in \gl_{\mu}(V)$ have Jordan-Chevalley decomposition
$L=S+N$ with $S,N \in \gl(V)$. Then $\gamma(L) = \mu L = \mu S + \mu N$
and so $\mu S + \mu N$ is the Jordan-Chevalley decomposition of
$\gamma(L)$. We also have $\gamma(L) = \gamma(S) + \gamma(N)$. Since
$\gamma(S) = \bfs^{-1}S\bfs$ or $\gamma(S) = \bfs^{-1}S^*\bfs$ we see that
$\gamma(S)$ is semi-simple. Furthermore $\gamma(N)^n = \gamma(N^n)=0$ and
so $\gamma(N)$ is nilpotent. Finally $[\gamma(S),\gamma(N)] = \sigma(\gamma)
\gamma([S,N]) = 0$. Thus $\gamma(S) + \gamma(N)$ is the
Jordan-Chevalley decomposition of $\gamma(L)$. Since this decomposition
is unique we have $\gamma(S)=\mu S$ and $\gamma(N) =\mu N$. \hfill $\Box$

\begin{remark}{\rm
Note that the subspaces of symmetric and skew-symmetric matrices in
$\gl(n,\fR)$, the $\pm 1$ eigenspaces of the anti-automorphism
$L \mapsto L^*$, both consist entirely of semi-simple matrices. Thus
for both these eigenspaces the Jordan-Chevalley decomposition is
trivial.\hfill $\rhd$}
\end{remark}

%
\subsubsection*{Coordinate Transformations}\label{sec:auttrafos}
We next look at coordinate transformations. If $L \in \gl(V)$ and $g \in
\Gl(V)$, then applying the coordinate change $g$ transforms $L$ to
$gLg^{-1}$. An automorphism $\phi_{\bfs}$ of $\gl(V)$ transforms to
$\phi_{g \bfs g^{-1}}$ and an anti-automorphism $\psi_{\bfs}$ to $\psi_{g
\bfs g^*}$. We therefore consider the action of $\Gl(V)$ on pairs $(L,\bfs)$
given by $(L,\bfs) \mapsto g . (L,\bfs)$, where $g . (L,\bfs)$ is a
shorthand for $(g L g^{-1}, g\bfs g^{-1})$ if $\gamma_{\bfs}$ is an
automorphism and $(g L g^{-1}, g \bfs g^*)$ if $\gamma_{\bfs}$ is an
anti-automorphism. Classifying pairs with respect to this action is
equivalent to first bringing $\bfs$ into normal form using any
transformation from $\Gl(V)$, and then classifying maps $L$ using only
transformations which preserve $\bfs$. The group of such transformations is
\begin{eqnarray*}
\Gl_{\bfs}^{+1} &=
    &\{g \in \Gl(V) \;|\; g\bfs g^{-1} = \bfs\} \\
\Gl_{\bfs}^{-1} &=
    &\{g \in \Gl(V) \;|\; g\bfs g^* = \bfs\}
\end{eqnarray*}
for automorphism and anti-automorphisms, respectively. The group $\Gl(V)$
can be indentified with the set of invertible elements in $\gl(V)$ and so
the (anti)-automorphism $\gamma_{\bfs}$ can also be regarded as an operator
on $\Gl(V)$. Thus $\gamma_{\bfs}(g) = \bfs^{-1}g \bfs$ for automorphisms
and $\gamma_{\bfs}(g) = \bfs^{-1} g^* \bfs$ for anti-automorphisms.  As
before we set $\sigma(\gamma_{\bfs})$ equal to $+1$ for automorphisms to $-1$
for anti-automorphims. Then the group of \emph{structure preserving}
transformations can be characterized as
\begin{displaymath}
\Gl_{\bfs}^{\sigma(\gamma_{\bfs})}(V) \ =\ \{g \in \Gl(V) \;|\;
     \gamma_{\bfs}(g) = g^{\sigma(\gamma_{\bfs})}\}.
\end{displaymath}
With help of the structure preserving transformation group we summarize the
discussion so far in the following lemma.
\begin{lemma}
\label{lem:orbit}
\textbf{\emph{(Orbit Lemma)}}\\
The $\Gl(V)$-orbit of the pair $(L,\bfs)$ is equivalent to
the $\Gl_{\bfs}^{\sigma(\gamma_{\bfs})}(V)$-orbit of $L$.
\end{lemma}
There is, however, a slightly larger transformation group that also
preserves the eigenspaces of $\gamma_{\bfs}$. This is the subgroup of
$\Gl(V)$ consisting of transformations which preserve the
(anti)-automorphism $\gamma = \gamma_{\bfs}$, rather than the structure map
$\bfs$ itself. We denote this group by
\begin{displaymath}
\Gl_{\gamma}(V) = \left\{ g \in \Gl(V)\;
\begin{array}{|lcll}
\phi_{g\bfs g^{-1}} & = & \phi_{\bfs}, & \mbox{if $\gamma = \phi_\bfs$ is an
automorphism}\\
\psi_{g\bfs g^*} & = & \psi_{\bfs}, & \mbox{if $\gamma = \psi_\bfs$ is an
anti-automorphism.}
\end{array}\right\}
\end{displaymath}
Equivalently, if $\phi_g(L) = g^{-1}Lg$ for $L\in \gl(V)$, then
\begin{equation} \label{eq:Glgammachar}
\Gl_{\gamma}(V) \ =\  \left\{g \in \Gl(V) \;|\;
\phi_g \circ \gamma = \gamma \circ \phi_g \right\}.
\end{equation}
In this paper we will use the groups $\Gl_{\bfs}^{\pm 1}(V)$ in the normal
form and unfolding theories. However in Lemma \ref{lem:int} we show that
the $\Gl_{\gamma}(V)$-orbit through $L \in \gl_\mu(V)$ is precisely the
intersection of the $\Gl(V)$-orbit in $\gl(V)$ with $\gl_\mu(V)$. Thus the
difference between the two groups is closely related to the phenomenon of
orbit splitting. The following proposition describes some of the elementary
properties of these transformation groups.
\begin{remark}{\rm
With a slight abuse of notation we define $\Sl(V) \defdas \{g \in \Gl(V)
\;|\; \det g =\pm 1\}$.  Note that the $\Gl(V)$ and $\Sl(V)$-orbits through
any $L$ are equal. We can therefore always work with the subgroups
$\Sl_{\bfs}^{\sigma(\gamma_{\bfs})}(V) = \Sl(V)\cap
\Gl_{\bfs}^{\sigma(\gamma_{\bfs})}(V)$ and $\Sl_{\gamma}(V) = \Sl(V)\cap
\Gl_{\gamma}(V)$ rather than $\Gl_{\bfs}^{\sigma(\gamma_{\bfs})}(V)$ and
$\Gl_{\gamma}(V)$ themselves.}
\end{remark}
\begin{proposition}
\label{pro:auttrafo}
Let $\gamma$ be an (anti)-automorphism of order two on $\gl(V)$ with
structure map $\bfs$. Let $\gl_{\mu}(V)$ be the eigenspace of $\gamma$ with
eigenvalue $\mu$.
\begin{enumerate}\abc
\item The groups $\Gl_{\bfs}^{\sigma(\gamma)}(V)$ and $\Gl_{\gamma}(V)$
preserve the eigenspace $\gl_{\mu}(V)$.
\item The Lie algebra of $\Gl_{\bfs}^{\sigma(\gamma)}(V)$ is the eigenspace
$\gl_{\sigma(\gamma)}(V)$.
\item The group $\Sl_{\bfs}^{\sigma(\gamma)}(V)$ is equal to either
the whole of $\Sl_{\gamma}(V)$ or to a normal subgroup of index two.
If $\dim V$ is odd
then $\Sl_{\bfs}^{\sigma(\gamma)}(V) = \Sl_{\gamma}(V)$.
If $\gamma$ is an automorphism then the same results hold
with
$\Gl_{\bfs}^{-1}(V)$ and  $\Gl_\gamma(V)$ in place of
$\Sl_{\bfs}^{-1}(V)$ and $\Sl_\gamma(V)$.

\end{enumerate}
\end{proposition}
Parts a) and b) of Proposition \ref{pro:auttrafo} also hold for
(anti)-automorphisms of any finite order.
\begin{example}
\label{ex:symmauts}{\rm
Let $\psi$ be the anti-automorphism $L \mapsto L^*$, for which the
structure map $\bfs$ is the identity map {\bf I}. Then
$\Gl_{\bfs}^{-1}(V) = \Sl_{\bfs}^{-1}(V) = \Sl_\gamma(V)$
is the group of orthogonal transformations. However $\Gl_\gamma(V)$
is the subgroup of $\Gl(V)$ consisting of elements $g$ such that
$gg^*$ is a scalar multiple of the identity. Thus part c) of
Proposition \ref{pro:auttrafo} does not hold with
$\Gl_{\bfs}^{-1}(V)$ and  $\Gl_\gamma(V)$ in place of
$\Sl_{\bfs}^{-1}(V)$ and $\Sl_\gamma(V)$.\hfill
$\rhd$}
\end{example}
\begin{example}{\rm
Consider the set of infinitesimally $R$-reversible maps on $\fR^2$ with
$R=\diag{1,-1}$ and let $\phi(L) = R^{-1}LR$. Then $\Sl_R^{+1}(V) =
\left\{\left(\begin{array}{cc}a&0\\0&b\end{array}\right) \;|\; ab = \pm
1\right\}$ has two cosets in $\Sl_{\gamma}(V)$. One coset is
$\Sl_R^{+1}(V)$, the other is $h\Sl_R^{+1}(V)$,
where $h=\left(\begin{array}{cc}0&1\\1&0\end{array}\right)$.
In this case it is also true that $\Gl_{\gamma}(V) =
\Gl_R^{+1}(V) \cup h\Gl_R^{+1}(V)$.
\hfill $\rhd$}
\end{example}
{\bf Proof of Proposition \ref{pro:auttrafo}.}
The proofs of parts \emph{a)} and \emph{b)} are straightforward calculations.

For part \emph{c)} we use the homomorphism $\rho : g \mapsto
\gamma(g)g^{-\sigma(\gamma)}$ of $\Gl(V)$. The characterization in equation
(\ref{eq:Glgammachar}) implies that for every element $g \in
\Gl_{\gamma}(V)$ the element $\rho(g) = \gamma(g)g^{-\sigma(\gamma)}$
commutes with every linear map $L \in \gl(V)$ and so must be a real nonzero
scalar multiple of the identity, $\gamma(g)g^{-\sigma(\gamma)} = \lambda
I$, say. Let $n = \dim V$. Then taking determinants, and noting that
$\det\gamma(g) = \det g$, implies that $\lambda^n = 1$ for automorphisms
and $\lambda^n = (\det g)^2$ for anti-automorphisms. It follows that
$\lambda = \pm 1$ for any $g \in \Gl_{\gamma}(V)$ if $\gamma$ is an
automorphisms, and for $g \in \Sl_{\gamma}(V)$ if $\gamma$ is an
anti-automorphism. Hence $\rho$ takes values in $\fZ_2 = \{\pm I\}$.  If
$n$ is odd then $\lambda = 1$ in both cases and $\rho$ is the constant
mapping to the identity. The results now follow.  \hfill $\Box$

%
\subsection{Reduction}\label{sec:autreduc}
We will now reduce to the semi-simple case using a method analogous
to that described in
Section \ref{sec:reducglv}. The next lemma follows almost immediately from
Lemma \ref{lem:reduc}.
\begin{lemma}
\label{lem:autreduc}
\textbf{\emph{(Reduction Lemma)}}\\ Let $L$ be a map in $\gl_{\mu}(V) = \{A
\in \gl(V) \;|\; \gamma_{\bfs}(A) = \mu A \}$, where $\gamma_{\bfs}$ is an
(anti)-automorphism of order two with structure map
$\bfs$. Furthermore let $L=S+N$ be the Jordan-Chevalley decomposition and
let $V_{\lambda}$ be an indecomposable $L$-invariant subspace. Then there
exists an indecomposable $\bfs$-invariant subspace $X_{\lambda} =
V_{\lambda} + \bfs V_{\lambda}$. Furthermore for each $X_{\lambda}$ there
exists an $S$-invariant complement $Y_{\lambda}$ of $NX_{\lambda}$ in
$X_{\lambda}$ such that $X_{\lambda}=Y_{\lambda} \oplus NX_{\lambda} =
Y_{\lambda} \oplus NY_{\lambda} \oplus \cdots \oplus
N^{n-1}Y_{\lambda}$. If $S$ is given on $Y_{\lambda}$, then $L$ is
determined on $X_{\lambda}$ up to similarity.
\end{lemma}
On $Y_{\lambda}$ we have a reduced (anti)-automorphism. For automorphisms
it is easy to see that the reduction of $\phi_{\bfs}$ is simply
$\phi_{\bfs}$ restricted to $Y_{\lambda}$ so that $Y_{\lambda} =
W_{\lambda} + \bfs W_{\lambda}$. Since $Y_{\lambda}$ is $\bfs$-invariant,
the normal form of $\bfs$ on $X_{\lambda}$ follows from the normal form of
$\bfs$ restricted to $Y_{\lambda}$.
\begin{remark}{\rm
Either $V_{\lambda} \cap \bfs V_{\lambda} = \{0\}$ or $V_{\lambda} \cap
\bfs V_{\lambda} = V_{\lambda}$. For example for infinitesimally
$R$-reversible maps we have $X_{\pm i \beta} = V_{\pm i
\beta}$, but $X_{\alpha} = V_{\alpha} \oplus R V_{\alpha} = V_{\alpha}
\oplus V_{-\alpha}$.\hfill $\rhd$}
\end{remark}

\commentaar{het volgende is in zekere zin uiterst onbevredigend!}
For anti-automorphisms the situation is somewhat more complicated. Recall
that an anti-automorphism comes from a non-degenerate bilinear form
$\omega$ which is either symmetric or anti-symmetric, $\omega(y,x) = \eps
\omega(x,y)$ with $\eps = \pm 1$. Then on $X_{\lambda}$ the map $L$ satisfies
$\omega(Lx,y) = \mu \omega(x,Ly)$. On $Y_{\lambda}$ we have a reduced
form.
\begin{lemma}
Let $L=S+N$ be the Jordan-Chevalley decomposition of $L \in \gl_{\mu}(V)$
on the indecomposable $\bfs$-invariant subspace $X_{\lambda}$. Let
$Y_{\lambda}$ be an $S$-invariant complement to $NX_{\lambda}$ in
$X_{\lambda}$. Then $\tau(x,y) = \omega(x,N^{n-1}y)$ is a non-degenerate
bilinear form on $Y_{\lambda}$ with $\tau(y,x) = \eps
\mu^{n-1}\tau(x,y)$. Furthermore $\tau(Sx,y) = \mu\tau(x,Sy)$.
\end{lemma}
See Appendix \ref{sec:gsbc} for a proof. The reduced anti-automorphism is
$\psi_{\bft}$ where the structure map $\bft$ is defined by $\tau(x,y) =
\inprod{x}{\bft y}_Y$ for all $x,y \in Y=Y_{\lambda}$. Then $Y_{\lambda} =
W_{\lambda} + \bft W_{\lambda}$. Here $Y_{\lambda}$ is \emph{not}
$\bfs$-invariant. However the freedom in choosing $Y_{\lambda}$ can be used
to put $\bfs$ into a  normal form with respect to the splitting
$X_{\lambda} =
Y_{\lambda} \oplus NY_{\lambda} \oplus \cdots \oplus N^{n-1}Y_{\lambda}$,
again see Appendix \ref{sec:gsbc}.

Now we obtain an unfolding lemma which uses the $\Gl(V)$-unfolding
from Lemma \ref{lem:unfo} as a starting point.
\begin{lemma}
\label{lem:autunfo}
\textbf{\emph{(Unfolding Lemma)}}\\ Let $\gamma$ be an (anti)-automorphism
and let $L \in \gl_{\mu}(V)$. Then the restriction of the $\Gl(V)$
centralizer unfolding of $L$ in $\gl(V)$ to $\gl_{\mu}(V)$
is equivalent to the
$\Gl_{\bfs}^{\sigma(\gamma)}(V)$ centralizer unfolding in $\gl_{\mu}(V)$.
\end{lemma}

\textbf{Proof.} Let $L \in \gl_{\mu}(V)$ and let $T_{\Gl(V)} = \{ UL-LU
\;|\; U \in \gl(V)\}$ be the tangent space at $L$ to the $\Gl(V)$-orbit of
$L$.  Let $N_{\Gl(V)}$ be its orthogonal complement with respect to an
inner product on $\gl(V)$. To simplify notation let $\G =
\Gl_{\bfs}^{\sigma(\gamma)}(V)$ and $\g = \gl_{\sigma(\gamma)}(V)$. Then
$T_{\G} = \{ UL-LU \;|\; U \in \g\}$ is the tangent space at $L$ to the
$\G$-orbit of $L$. Let $N_{\G}$ be its orthogonal complement in
$\gl_{\mu}(V)$. By shifting $L$ to the origin in $\gl(V)$ we have $\gl(V) =
T_{\Gl(V)} \oplus^{\perp} N_{\Gl(V)}$ and $\gl_{\mu}(V) = T_{\G}
\oplus^{\perp} N_{\G}$. Now let $\Pi$ be the orthogonal projection onto
$\gl_{\mu}(V)$. Then $\Pi(T_{\Gl(V)}) = \{UL-LU \;|\; U \in \gl(V),
\gamma(UL-LU) = \mu(UL-LU)\} = \{ UL-LU \;|\; U \in \g\} = T_{\G}$, since
$\gamma(UL-LU) = \mu(UL-LU)$ iff $U \in \g$, and
\begin{displaymath}
\gl_{\mu}(V) = \Pi(\gl(V)) = \Pi(T_{\Gl(V)} \oplus^{\perp} N_{\Gl(V)}) =
\Pi(T_{\Gl(V)}) \oplus^{\perp} \Pi(N_{\Gl(V)}) = T_{\G} \oplus^{\perp}
\Pi(N_{\Gl(V)}).
\end{displaymath}
Since orthogonal complements are unique, we have $\Pi(N_{\Gl(V)}) =
N_{\G}$.\hfill $\Box$

%
\subsection{Orbit Splitting}\label{sec:osplit}
The transformation group $\Gl_{\gamma}(V)$ preserves $\gl_{\mu}(V)$ as a
linear space and so maps $\Gl_{\bfs}^{\sigma(\gamma)}(V)$-orbits into
$\Gl_{\bfs}^{\sigma(\gamma)}(V)$-orbits.  The difference between
the groups $\Gl_{\gamma}(V)$ and $\Gl_{\bfs}^{\sigma(\gamma)}(V)$
gives rise to \emph{splitting} of orbits. This
means that a $\Gl(V)$-orbit in $\gl(V)$ splits into different
$\Gl_{\bfs}^{\sigma(\gamma)}(V)$-orbits when intersected with $\gl_\mu(V)$.
In particular orbit splitting
implies that eigenvalues and Jordan structure no longer suffice to
characterize $\Gl_{\bfs}^{\sigma(\gamma)}(V)$-orbits. Examples of additional
invariants are \emph{symplectic} and \emph{reversible signs}
\cite{bc3,hov}.
\begin{lemma}
\label{lem:int}
\textbf{\emph{(Intersection Lemma)}}\\ Let $V$ be an indecomposable
$\bfs$-invariant space and $\gamma_{\bfs}$ an (anti)-automorphism of
order two with structure map $\bfs$. Let $L \in \gl_\mu(V)$. Then
\begin{displaymath}
\Orb_{\Gl(V)}(L) \ \cap\ \gl_\mu(V)\ = \ \Orb_{\Gl_\gamma(V)}(L).
\end{displaymath}
\end{lemma}

\textbf{Proof.} If $M \in \Orb_{\Gl_\gamma(V)}(L)$ then clearly $M \in
\Orb_{\Gl(V)}(L) \cap \gl_\mu(V)$. Next we assume $M \in \Orb_{\Gl(V)}(L)
\cap \gl_\mu(V)$. The $\Gl(V)$-orbits of $L$ and $M$ are, of course, the
same. From the Orbit Lemma \ref{lem:orbit} we know that the
$\Gl^{\sigma(\gamma)}_\bfs(V)$-orbits of $L$ and $M$ correspond to the
$\Gl(V)$-orbits of the pairs $(L,\bfs)$ and $(M,\bfs)$. These orbits are
not necessarily the same. From the $\Gl(V)$-classification of pairs in
Section \ref{sec:nfunfo} we see that the normal form of $(L,\bfs)$ is
either $(L_0,\bfs_0)$ or $(L_0,\eps\bfs_0)$ with $\eps=\pm 1$. That is
either eigenvalues and Jordan structure determine the orbit or there is an
additional sign. In that case $(L_0,\bfs_0)$ and $(L_0, -\bfs_0)$ are
\emph{not} equivalent, there exists no $g \in \Gl(V)$ such that $g
. (L_0,\bfs_0) = (L_0,-\bfs_0)$. The normal form of $(M,\bfs)$ is
$(L_0,\eps\bfs_0)$ and in the absence of a sign we set $\eps=1$. If there
is no sign then there are $g_1, g_2 \in \Gl(V)$ such that $g_1 . (L,\bfs) =
(L_0,\bfs_0)$ and $g_2 . (M,\bfs) = (M_0,\bfs_0)$. But then we have
$(M,\bfs) = g_2^{-1}g_1 . (L,\bfs)$ and thus $g_2^{-1}g_1 \in
\Gl_\bfs^{\sigma(\gamma)}(V) \subset \Gl_{\gamma}(V)$. If there is a sign
we perform the same computation to find $g_1, g_2 \in \Gl(V)$ such that
$(M,\bfs) = g_2^{-1}g_1 . (L,\eps\bfs)$ and thus $g_2^{-1}g_1 \in
\Gl_{\gamma}(V)$. Hence we may conclude that $M \in
\Orb_{\Gl_\gamma(V)}(L)$. \hfill $\Box$

The following result is an immediate corollary of Proposition
\ref{pro:auttrafo} and Lemma \ref{lem:int}.
\begin{theorem}
\label{the:split}
\textbf{\emph{(Orbit Splitting Theorem)}}\\ If $L \in \gl_{\mu}(V)$ then the
$\Gl(V)$-orbit of $L$ in $\gl(V)$ intersects $\gl_{\mu}(V)$ in at most two
$\Gl_{\bfs}^{\sigma(\gamma)}(V)$-orbits.
\end{theorem}
\begin{example}{\rm \label{ex:rev}
Recall from Example \ref{ex:motiv1} that the space of infinitesimally
$R$-reversible maps is defined by $\gl_{-R}(V) \defdas \{L \in \gl(V)\;|\;
\phi_R(L)=-L\}$, where $\phi_R$ is an automorphism defined by $\phi_R(L)
\defdas R^{-1}LR$ with \emph{structure map} $R$ satisfying $R^2=I$, $R \neq
\pm I$. The structure map $R$ has eigenvalues $+1$ and $-1$ with
corresponding eigenspaces $E_+=\{a \in V \;|\; Ra = a\}$ and $E_-=\{b \in V
\;|\; Rb =-b\}$. For an infinitesimally $R$-reversible map $L$ we have $L
E_+ \subseteq E_-$ and $L E_- \subseteq E_+$. Let $V$ be an indecomposable
$L$,$R$-invariant subspace. Then eigenvectors of $L$ for eigenvalue zero
either belong to $E_+$ or to $E_-$. The new invariant, called the
\emph{reversible sign} indicates to which eigenspace such an eigenvector
belongs.\hfill $\rhd$}
\end{example}
\begin{example}{\rm \label{ex:sympl}
Infinitesimally symplectic linear maps are defined by $\sp(V) \defdas \{L
\in \gl(V) \;|\; \omega(Lx,y) + \omega(x,Ly) = 0,\; \forall x,y \in V\}$,
where $\omega$ is a non-degenerate skew symmetric bilinear form on $V$, see
Example \ref{ex:motiv2}.  With help of an inner product
$\inprod{\cdot}{\cdot}$ on $V$ we can find a map $J$ satisfying $J^*=-J$
and $J^2=-I$ such that $\omega(x,y) = \inprod{x}{Jy}$. Then $\sp(V) = \{L
\in \gl(V) \;|\; \psi_J(L) = -L\}$, where $\psi_J$ is the anti-automorphism
defined by $\psi_J(L) \defdas J^{-1}L^*J$. Here $J$ is the \emph{structure
map}. Since $\gl_{\sigma(\psi_J)}(V) = \gl_{-1}(V) = \sp(V)$, Proposition
\ref{pro:auttrafo} says that $\sp(V)$ is a Lie algebra, as we already
knew. The invertible transformations that preserve the symplectic structure
are exactly those $g$ for which $\omega(gx,gy) = \omega(x,y)$. Indeed
according to Proposition \ref{pro:auttrafo} they are given by $\Gl_1(V) =
\{g \in \Gl(V) \;|\; \psi_J(g) = g^{-1}\}$, which is equivalent to
$g^*Jg=J$, reflecting the transformation rule for bilinear forms. The order
of $\psi_J$ is two. So if splitting of an orbit occurs there are at most
two inequivalent normal forms. A well known example is the distinction
between 1:1 resonance and 1:-1 resonance. In both cases there are double
eigenvalues $\pm \bfi$. But there is an additional invariant, namely a
\emph{symplectic} sign distinguishing the two cases. These signs are
intimately related to the Morse index of the corresponding quadratic
Hamilton functions. In Hamiltonian systems a single pair of complex
conjugate imaginary eigenvalues is forced to remain on the imaginary axis
when parameters of the system vary. When two such pairs meet they may
remain on the imaginary axis, which is called \emph{passing} or they may
move into the complex plane which is called \emph{splitting}. Computing the
unfoldings in the 1:1 case one finds passing of imaginary eigenvalues (see
Cotter \cite{cot}, Galin \cite{gal} and Ko\c{c}ak \cite{koc}), which is a
codimension three phenomenon, see table \ref{tab:unfo5678} type 8d. In the
1:-1 case one finds splitting of imaginary eigenvalues (see van der Meer
\cite{mee}), which is a codimension one phenomenon, see table
\ref{tab:unfo5678} type 8c.\hfill $\rhd$}
\end{example}

%
\section{Normal Forms and Examples of Unfoldings}\label{sec:nfunfo}
%

%
\subsection{Normal Forms}\label{sec:nfs}
In the present setting, where $\gamma_{\bfs}$ is an (anti)-automorphism of
order two with structure map $\bfs$, we can classify maps satisfying
$\gamma_{\bfs}(L) = \mu L$ into the eight different types listed
in Table \ref{tab:types}.
\begin{table}[hbtp]
\begin{displaymath}
\begin{array}{|c|c|c|c|l|}\hline&&&&\\[-2ex]
\mbox{type} & \gamma_{\bfs} &  \bfs^*   & \mu    & L \\[1ex]\hline&&&&\\[-2ex]
1&\phi_{\bfs}   & \phm \bfs & \phm 1 & \mbox{$\bfs$-equivariant of type
$\fR$}\\[1.5ex]
2&\phi_{\bfs}   & \phm \bfs &     -1 & \mbox{$\bfs$-reversible of type
$\fR$}\\[1.5ex]
3&\phi_{\bfs}   &     -\bfs & \phm 1 & \mbox{$\bfs$-equivariant of type $\fC$
or $\fC$-linear}\\[1.5ex]
4&\phi_{\bfs}   &     -\bfs &     -1 & \mbox{$\bfs$-reversible of type $\fC$
or $\fC$-semi-linear}\\[1.5ex]\hline&&&&\\[-2ex]
5&\psi_{\bfs}   & \phm \bfs & \phm 1 & \mbox{symmetric}\\[1.5ex]
6&\psi_{\bfs}   & \phm \bfs &     -1 & \mbox{anti-symmetric}\\[1.5ex]
7&\psi_{\bfs}   &     -\bfs & \phm 1 & \mbox{anti-symplectic}\\[1.5ex]
8&\psi_{\bfs}   &     -\bfs &     -1 & \mbox{symplectic}\\[1.5ex]\hline
\end{array}
\end{displaymath}
\caption{\textit{Eight types of maps satisfying $\gamma_{\bfs}(L) = \mu L$,
when $\gamma_{\bfs}$ is an (anti)-automorphism of order
two.\label{tab:types}}}
\end{table}

Applying the Reduction Lemma \ref{lem:autreduc} it is straightforward to
obtain normal forms for semi-simple maps of the types listed in Table
\ref{tab:types}. If $\gamma_{\bfs}$ is an anti-automorphism we apply the
construction of Appendix \ref{sec:gsbc} to put the structure map $\bfs$
into normal form. Finally we reduce to the smaller space $Y$ on which we
consider the semi-simple part $S$ of $L$ and the reduced structure map
$\bft$ (see Appendix \ref{sec:gsbc}), such that $S \in \gl_{\mu}(Y) = \{A
\in \gl(Y) \;|\;\gamma_{\bft}(A)=\mu A\}$. Since we consider semi-simple
maps in their own right we denote the reduced structure map again by
$\bfs$. Note that in view of the Orbit Lemma \ref{lem:orbit} the normal
forms in Tables \ref{tab:nfs1234} and \ref{tab:nfs5678} can be regarded as
representatives of $\Gl^{\sigma(\gamma)}_\bfs(Y)$-orbits of $S$ once a
choice for $\bfs$ has been made, but they can also be seen as
representatives of $\Gl(Y)$-orbits of the pair $(S,\bfs)$.
\begin{theorem}
\label{the:nfs}
\textbf{\emph{(Normal Form Theorem)}}\\ Let $\gamma_{\bfs}$ be an
(anti)-automorphism of order two with structure map $\bfs$ and let $L \in
\gl_{\mu}(X)$ be a linear map in one of the eigenspaces of $\gamma_{\bfs}$,
where $X$ is an indecomposable $(L,\bfs)$-invariant space. Furthermore let
$L=S+N$ be the Jordan-Chevalley decomposition of $L$. Then the normal form
of $L$ on $X$ is determined by the normal form of $S$ on $Y$, where $Y
\subset X$ is a reduced space as in Proposition \ref{lem:autreduc}. The
normal forms of $S$ are listed in Tables \ref{tab:nfs1234} and
\ref{tab:nfs5678}.
\end{theorem}
\begin{remarks}{\rm \label{rem:signs}
Let us make some remarks on the various cases in Tables \ref{tab:nfs1234}
and \ref{tab:nfs5678}.
\begin{enumerate}\parskip 0pt \itemsep 0pt
\item
There are only two cases where there are no signs at all, namely
infinitesimally $\bfs$-reversible semi-simple maps of type $\fC$ and
anti-symplectic semi-simple maps.
\item In each of the other cases the sign has a geometric meaning. Let $S$
be a semi-simple map in an eigenspace of an automorphism, see Table
\ref{tab:types}.
\begin{enumerate}\abc\parskip 0pt \itemsep 0pt
\item For $\bfs$-equivariant maps of type $\fR$ there are two isotypic
components, labelled by the eigenvalues of $\bfs$. The sign indicates on
which of the two isotypic components $L$ acts. Similarly for
$\bfs$-equivariant maps of type $\fC$. Here the sign disappears for real
eigenvalues of $L$. The latter are forced to be double and are non-generic.
\item For infinitesimally $\bfs$-reversible maps of type $\fR$ the sign for
zero eigenvalues of $L$ indicates to which eigenspace of $\bfs$ the
corresponding eigenvector belongs.
\item In case of symmetric and anti-symmetric maps the sign is related to
the signature of the underlying bilinear form. If the signature is denoted
by $(p,q)$, that is
\begin{displaymath}
\inprod{x}{y} = \sum_{i=1}^p x_i \, y_i - \sum_{i=1}^{q} x_{p+i} \, y_{p+i},
\end{displaymath}
then for signature $(n,0)$ or $(0,n)$, symmetric maps only have real
eigenvalues whereas anti-symmetric maps only have imaginary eigenvalues. In
this case the sign distinguishes between $(n,0)$ and $(0,n)$.
\item For infinitesimally symplectic maps the structure map itself is
symplectic. The dynamical interpretation is that the structure map defines
a preferred direction of rotation. This gives a sign for imaginary
eigenvalues of an infinitesimally symplectic map.
\end{enumerate}
\item In case of types 3 and 4 the structure map $\bfs$ satisfies
$\bfs^*=-\bfs$ and $\bfs^2=-I$. Therefore $\bfs$ defines a complex
structure on $X$. Maps of type 3 commute with $\bfs$, but maps of type 4
anti-commute with $\bfs$. Consequently maps of type 3 can be considered as
$\fC$-linear maps whereas maps of type 4 can be considered as
$\fC$-semi-linear maps.
\item In case of anti-automorphisms the reduced structure map $\bft$ can be
either symmetric or skew in eigenspaces with $\mu=-1$. This means that a
reduced semi-simple infinitesimally symplectic map is either
infinitesimally symplectic or anti-symmetric with respect to the reduced
structure map. Similarly a reduced semi-simple anti-symmetric map is either
anti-symmetric or infinitesimally symplectic.\hfill $\rhd$
\end{enumerate}
}
\end{remarks}
\commentaar{Repercussies voor de signatuur}

\begin{table}[htb]
\begin{displaymath}
\begin{array}{|c|c|c|c|c|c|}\hline&&&&&\\[-2ex]
\mbox{type} & \lambda & Y_{\lambda} & \bfs & S & \mbox{remarks}
\\[1ex]\hline&&&&&\\[-2ex]
1 & \alpha  & W_{\alpha}  & \pm 1 & \alpha  & \mbox{sign}\\[1.5ex]
1 & \alpha \pm i\beta & W_{\alpha \pm i\beta} &  \pm I_2 &
\left(\begin{array}{rr}\alpha&-\beta\\\beta&\alpha\end{array}\right) &
\mbox{sign, $\beta > 0$}\\[2ex]\hline&&&&&\\[-2ex]
2 & 0 & W_0 & \pm 1 & 0 &\mbox{sign}\\[1.5ex]
2 & \pm \alpha & W_{\alpha} \oplus W_{-\alpha} & T &
\left(\begin{array}{rr}\alpha&0\\0&-\alpha\end{array}\right) &
\alpha > 0\\[3ex]
2 & \pm i\beta & W_{\pm i\beta} & R &
\left(\begin{array}{rr}0&-\beta\\\beta&0\end{array}\right) &
\beta > 0\\[3ex]
2 & \pm \alpha\pm i\beta & W_{\alpha\pm i\beta} \oplus W_{-\alpha\pm i\beta} &
\left(\begin{array}{rr}0&I_2\\I_2&0\end{array}\right) &
\left(\begin{array}{rrrr}\alpha&-\beta&&\\\beta&\alpha&&\\
&&-\alpha&\beta\\&&-\beta&-\alpha\end{array}\right) &
\alpha > 0,\;\beta > 0\\[5ex]\hline&&&&&\\[-2ex]
\commentaar{
3 & \alpha\alpha & W_{\alpha} \oplus \bfs W_{\alpha} & J &
\left(\begin{array}{rr}\alpha&0\\0&\alpha\end{array}\right) &\\[3ex]
}
3 & \alpha \pm i\beta & W_{\alpha \pm i\beta} & \pm J &
\left(\begin{array}{rr}\alpha&-\beta\\\beta&\alpha\end{array}\right) &
\mbox{sign, $\beta \geq 0$}\\[2ex]\hline&&&&&\\[-2ex]
4 & \pm \alpha & W_{\alpha} \oplus W_{-\alpha} & J &
\left(\begin{array}{rr}\alpha&0\\0&-\alpha\end{array}\right) &
\alpha \geq 0\\[3ex]
4 & \pm \alpha\pm i\beta & W_{\alpha\pm i\beta} \oplus W_{-\alpha\pm i\beta} &
\left(\begin{array}{rr}0&-I_2\\I_2&0\end{array}\right) &
\left(\begin{array}{rrrr}\alpha&-\beta&&\\\beta&\alpha&&\\
&&-\alpha&\beta\\&&-\beta&-\alpha\end{array}\right) &\alpha \geq 0,\;
\beta > 0\\[5ex]\hline
\end{array}
\end{displaymath}
\caption{\textit{Normal forms for semi-simple maps $S$ in the eigenspace of
an automorphism $\phi_{\bfs}$. The type refers to the types in Table
\ref{tab:types}. $W_{\lambda}$ is an indecomposable $S$-invariant space on
which $S$ has eigenvalue $\lambda$ in the sense of Theorem
\ref{the:vlambda} and $Y_{\lambda} = W_{\lambda} + \bfs W_{\lambda}$ is an
indecomposable $S,\bfs$-invariant space. In the characterization of the
structure map $\bfs$ we the standard matrices $I$, $I_2$, $R$, $T$ and $J$
as in the proof of Theorem \ref{the:nfs}. The \emph{signs} in the tables
indicate that the normal forms for $+1$ and $-1$ are inequivalent. Since in
each case there are at most two possibilities this shows that orbit
splitting occurs as stated in Theorem \ref{the:split}. The relevance of the
signs becomes most obvious when constructing unfoldings. We put the sign in
the structure map $\bfs$, but for nonzero eigenvalues we can also put the
sign in the map $S$.}
\label{tab:nfs1234}}
\end{table}

\begin{table}[htb]
\begin{displaymath}
\begin{array}{|c|c|c|c|c|c|}\hline&&&&&\\[-2ex]
\mbox{type} & \lambda & Y_{\lambda} & \bfs & S & \mbox{remarks}
\\[1ex]\hline&&&&&\\[-2ex]
5 & \alpha & W_{\alpha} & \pm 1 & \alpha &\mbox{sign}\\[1.5ex]
5 & \alpha \pm i\beta & W_{\alpha \pm i\beta} & R &
\left(\begin{array}{rr}\alpha&-\beta\\\beta&\alpha\end{array}\right) &
\beta > 0\\[2ex]\hline&&&&&\\[-2ex]
6 & 0 & W_0 & \pm 1 & 0 & \mbox{sign}\\[1.5ex]
6 & \pm \alpha & W_{\alpha} \oplus W_{-\alpha} & T &
\left(\begin{array}{rr}\alpha&0\\0&-\alpha\end{array}\right) &
\alpha > 0\\[3ex]
6 & \pm i\beta & W_{\pm i\beta} &
 \pm I_2 & \left(\begin{array}{rr}0&-\beta\\\beta&0\end{array}\right) &
\mbox{sign, $\beta>0$}\\[3ex]
6 & \pm \alpha\pm i\beta & W_{\alpha\pm i\beta} \oplus W_{-\alpha\pm i\beta} &
\left(\begin{array}{rr}0&I_2\\I_2&0\end{array}\right) &
\left(\begin{array}{rrrr}\alpha&-\beta&&\\\beta&\alpha&&\\
&&-\alpha&-\beta\\&&\beta&-\alpha\end{array}\right) & \alpha > 0,\;
\beta > 0\\[5ex]\hline&&&&&\\[-2ex]
7 & \alpha\alpha & W_{\alpha} \oplus \bfs W_{\alpha} & J &
\left(\begin{array}{rr}\alpha&0\\0&\alpha\end{array}\right) &\\[3ex]
7 & \alpha \pm i\beta & W_{\alpha\pm i\beta} \oplus \bfs W_{\alpha\pm i\beta} &
\left(\begin{array}{rr}0&-I_2\\I_2&0\end{array}\right) &
\left(\begin{array}{rrrr}\alpha&-\beta&&\\\beta&\alpha&&\\
&&\alpha&\beta\\&&-\beta&\alpha\end{array}\right) & \beta > 0
\\[5ex]\hline&&&&&\\[-2ex]
8 & \pm \alpha & W_{\alpha} \oplus W_{-\alpha} & J &
\left(\begin{array}{rr}\alpha&0\\0&-\alpha\end{array}\right) &
\alpha \geq 0\\[3ex]
8 & \pm i\beta & W_{\pm i\beta} & \pm J &
\left(\begin{array}{rr}0&-\beta\\\beta&0\end{array}\right) &
\mbox{sign, $\beta>0$}\\[3ex]
8 & \pm \alpha\pm i\beta & W_{\alpha\pm i\beta} \oplus W_{-\alpha\pm i\beta} &
\left(\begin{array}{rr}0&-I_2\\I_2&0\end{array}\right) &
\left(\begin{array}{rrrr}\alpha&-\beta&&\\\beta&\alpha&&\\
&&-\alpha&-\beta\\&&\beta&-\alpha\end{array}\right) & \alpha > 0,\;
\beta >0\\[5ex]\hline
\end{array}
\end{displaymath}
\caption{\textit{Normal forms for semi-simple maps in eigenspaces of
anti-automorphisms. See the caption of Table \ref{tab:nfs1234} for an
explanation of the notation.\label{tab:nfs5678}}}
\end{table}

We conclude this section with a proof of the Normal Form Theorem.

\textbf{Proof of Theorem \ref{the:nfs}.} We will not give proofs for all
the different cases for they are very similar to each other. Instead we
give some representative proofs for automorphisms and anti-automorphisms
with and without signs. The types refer to Table \ref{tab:types}. We use
the following standard matrices.
\begin{displaymath}
I=I_2=\left(\begin{array}{rr}1&0\\0&1\end{array}\right),\;
R=\left(\begin{array}{rr}1&0\\0&-1\end{array}\right),\;
T=\left(\begin{array}{rr}0&1\\1&0\end{array}\right),\;
J=\left(\begin{array}{rr}0&-1\\1&0\end{array}\right).
\end{displaymath}
\begin{enumerate}\parskip 0pt \itemsep 0pt
\item \underline{Type 3; complex eigenvalues; $\bfs^2=-I$.} The maps $S$
and $\bfs$ commute, both have complex eigenvalues, and so $Y_{\alpha \pm
i\beta} = W_{\alpha \pm i\beta}$. Let $e$ be any vector in $Y_{\alpha \pm
i\beta}$ and define $f \defdas \frac{1}{\beta}(S-\alpha)e$.  Then
$\basis{e,f}$ is a basis of $Y_{\alpha \pm i\beta}$ on which we have $S =
\alpha I + \beta J$ and $\bfs = J$. Note that we can change the sign of
$\beta$ in $S$ by applying the transformation $R$ or $T$. Then the sign of
$\bfs$ changes as well. Therefore we assume $\beta > 0$ and put the sign in
$\bfs = \pm J$.
\item \underline{Type 4; real eigenvalues; $\bfs^2=-I$.} The maps $S$ and
$\bfs$ anti-commute.  If $e \in W_{\alpha}$ then $\bfs e \in W_{-\alpha}$,
so $Y_{\alpha} = W_{\alpha} \oplus W_{-\alpha}$. Then $\basis{e,\bfs e}$ is
a basis of $Y_{\alpha}$ and we have $S = \alpha R$ and $\bfs = J$. Note
that $\bfs S \bfs^{-1} = -S$, so we may assume that $\alpha \geq 0$. There
is no sign in this case.
\item \underline{Type 5; complex eigenvalues; $\bfs^2=I$.} Here $Y_{\alpha
\pm i\beta} = W_{\alpha \pm i\beta}$. Let $e$ be any vector in $Y_{\alpha
\pm i\beta}$ and define $f \defdas \frac{1}{\beta}(S-\alpha)e$. Then
$\basis{e,f}$ is a basis of $Y_{\alpha \pm i\beta}$. On this basis $S =
\alpha I + \beta J$. Let $\tau$ be the reduced form on $Y_{\alpha \pm
i\beta}$, that is $\tau(x,y)=\inprod{x}{\bfs y}$. Then indeed $\tau(e,f) =
\frac{1}{\beta}(e,Se) = \frac{1}{\beta}(Se,e) = \tau(f,e)$. If $u$ is a
vector in $Y_{\alpha \pm i\beta}$ such that $\bfs u = u$, then $\bfs
(S-\alpha)u = -(S-\alpha) \bfs u = -(S-\alpha) u$. So $\bfs$ is indeed
symmetric and has eigenvalues $\pm 1$ on $Y_{\alpha \pm i\beta}$. Thus
there exists an orthogonal transformation such that
$\basis{u,\frac{1}{\beta}(S-\alpha)u}$ is a new basis of $Y_{\alpha \pm
i\beta}$ and $S = \alpha I + \beta J$, $\bfs = R$. Note that $\bfs S
\bfs^{-1} = \alpha I - \beta J$, so we may assume $\beta \geq 0$ and there
is no sign.
\item \underline{Type 6; imaginary eigenvalues; $\bfs^2=I$.} Again
$Y_{\alpha \pm i\beta} = W_{\alpha \pm i\beta}$. Let $e$ be any vector in
$Y_{\pm i\beta}$ and define $f \defdas \frac{1}{\beta}Se$.  Then
$\basis{e,f}$ is a basis of $Y_{\pm i\beta}$. On this basis $S = \beta
J$. The reduced form $\tau$ on $Y_{\alpha \pm i\beta}$ is symmetric, so
$\tau(e,f)=\tau(f,e)$. On the other hand $\tau(e,f) =
\frac{1}{\beta}\tau(e,Se) = -\frac{1}{\beta}\tau(Se,e) = -\tau(f,e)$, so
$\tau(e,f)=0$. Furthermore $\tau(f,f) = \frac{1}{\beta^2}\tau(Se,Se) =
-\frac{1}{\beta^2}\tau(S^2e,e) = \tau(e,e)$. Thus on the basis
$\basis{e,f}$ we have $S = \beta J$, $\bfs = \pm I$ and we may assume that
$\beta > 0$.\hfill $\Box$
\end{enumerate}

%
\subsection{Examples of Unfoldings}\label{sec:unfos}
Our aim is to present a list of low codimension unfoldings. We will
concentrate on unfoldings of zero eigenvalues. Since such unfoldings give a
parametrization of the full space of maps we also get information on the
unfoldings of real, imaginary and complex eigenvalues. The classification
is complete up to codimension two. Where appropriate we include some higher
codimension cases to show the consequences of signs. In constructing
unfoldings we use Lemmas \ref{lem:autunfo} and \ref{lem:unfo} and we summarize
the results in a theorem.
\begin{theorem}
\label{the:unfo}
\textbf{\emph{(Unfolding Theorem)}}\\ Let $\gamma_{\bfs}$ be an
(anti)-automorphism of order two with structure map $\bfs$ and let $L \in
\gl_{\mu}(X)$ be a linear map, with $X$ an indecomposable
$(L,\bfs)$-invariant space. Furthermore let $L=S+N$ be the Jordan-Chevalley
decomposition of $L$. Then the unfoldings of $L$ up to codimension two are
listed in Tables \ref{tab:unfo1234} and \ref{tab:unfo5678}.
\end{theorem}

Here we consider real maps with real parameters. In such maps simple real
eigenvalues are forced to remain on the real axis when parameters are
varied. Maps which are elements of an eigenspace corresponding to
eigenvalue $-1$ of an (anti)-automorphism have eigenvalues which come in
complex conjugate $(\lambda,\bar{\lambda})$ \emph{and} opposite
$(\lambda,-\lambda)$ pairs. Thus a pair of simple opposite imaginary
eigenvalues is forced to remain on the imaginary axis when parameters are
varied. However there may also be collisions of such eigenvalues on the
real or imaginary axis as the parameter varies. Generically eigenvalues
\emph{split} into the complex plane at collisions, but in the presence of
signs generic \emph{passing} also occurs.

\begin{table}[htb]
\begin{displaymath}
\begin{array}{|c|c|c|c|c|c|}\hline&&&&&\\[-2ex]
\mbox{type} & \lambda & X_{\lambda} & \bfs & L(\nu) & \mbox{codim}
\\[1ex]\hline&&&&&\\[-2ex]
1a & 0^+ & \basis{a} & 1 & \nu & 1\\[1.5ex]
1b & (0^+)^2 & \basis{a,Na} & I_2 &
\left(\begin{array}{cc}\nu_1&\nu_2\\1&\nu_1\end{array}\right) & 2\\[3ex]
1c & (0^+)(0^-) & \basis{a,b} & R &
\left(\begin{array}{cc}\nu_1&0\\0&\nu_2\end{array}\right) & 2\\[3ex]
1d & (0^+)(0^+) & \basis{a_1,a_2} & I_2 &
\left(\begin{array}{cc}\nu_1&\nu_2\\\nu_3&\nu_4\end{array}\right) & 4
\\[2ex]\hline&&&&&\\[-2ex]
2a & 0^+ & \basis{a} & 1 & 0 & 0\\[1.5ex]
2b & (0^+)^2 & \basis{Nb,b} & R &
\left(\begin{array}{cc}0&1\\\nu&0\end{array}\right) & 1\\[3ex]
2c & (0^+)(0^-) & \basis{a,b} & R &
\left(\begin{array}{cc}0&\nu_1\\\nu_2&0\end{array}\right) & 2\\[3ex]
2d & (0^+)^3 & \basis{a,N^2a,Na} &
\left(\begin{array}{r|r}I_2&\\\hline&-1\end{array}\right) &
\left(\begin{array}{rr|r}&&\nu\\&&1\\\hline1&\nu&\end{array}\right) & 1\\[4ex]
2e & (0^+)^2(0^+) & \basis{a,Nb,b} &
\left(\begin{array}{r|r}I_2&\\\hline&-1\end{array}\right) &
\left(\begin{array}{rr|r}&&0\\&&1\\\hline\nu_2&\nu_1&\end{array}\right) &
2\\[4ex]
2f & (0^+)^2(0^-) & \basis{Nb_1,b_1,b_2} &
\left(\begin{array}{r|r}1&\\\hline&-I_2\end{array}\right) &
\left(\begin{array}{r|rr}&1&0\\\hline\nu_1&&\\\nu_2&&\end{array}\right) &
2\\[4ex]\hline&&&&&\\[-2ex]
3a & 00 & \basis{e, f} & J &
\left(\begin{array}{rr}\nu_1&-\nu_2\\\nu_2&\nu_1\end{array}\right) & 2\\[3ex]
3b & (\pm i\beta)^2 & \basis{e, f, Ne, N f} &
\left(\begin{array}{cc}J&\\&J\end{array}\right) &
\left(\begin{array}{rrrr}
\nu_1&-\beta-\nu_2&\nu_3&-\nu_4\\
\beta+\nu_2&\nu_1&\nu_4&\nu_3\\
1&0&\nu_1&-\beta-\nu_2\\
0&1&\beta+\nu_2&\nu_1
\end{array}\right) & 4\\[5ex]
3c & (\pm i\beta)(\pm i\beta) & \basis{e, \bfs e, f, -\bfs f} &
\left(\begin{array}{cc}J&\\&-J\end{array}\right) &
\left(\begin{array}{rrrr}
\nu_1&-\beta-\nu_2&&\\
\beta+\nu_2&\nu_1&&\\
&&\nu_3&-\beta-\nu_4\\
&&\beta+\nu_4&\nu_3
\end{array}\right) & 4\\[5ex]\hline&&&&&\\[-2ex]
4 & 0 & \basis{e, f} & J &
\left(\begin{array}{rr}\nu_1&\nu_2\\\nu_2&-\nu_1\end{array}\right) & 2
\\[2ex]\hline
\end{array}
\end{displaymath}
\caption{\textit{Unfoldings in eigenspaces of automorphisms. The notation
is similar to Tables \ref{tab:nfs1234} and \ref{tab:nfs5678}. Again the
type refers to the types in Table \ref{tab:types}. Here $X_{\lambda}$ is
the $(L,\bfs)$-invariant space on which $L$ has eigenvalue $\lambda$ in the
sense of Theorem \ref{the:vlambda} and Lemma \ref{lem:recons}. Eigenvalues
are denoted $\lambda^n$ when their multiplicity is $n$. We use brackets
when signs are present. For example $(0^+)^2(0^-)$ in type 2f means three
eigenvalues zero, one with muliplicity 2 and sign $+1$ and one with
multiplicity 1 and sign $-1$. Basis vectors in the tables are such that
$\bfs a = a$ and $\bfs b =-b$ for types 1 and 2. In all other cases $e$ and
$f$ are vectors in the complement of $NX_{\lambda}$ in $X_{\lambda}$ where
$f$ is generated by the semi-simple part $S$ of $L$ or the structure map
$\bfs$. For example $f = \frac{1}{\beta}Se$ in type 3b, but $f = \bfs e$ in
type 3a.\label{tab:unfo1234}} }
\end{table}

\begin{table}[htb]
\begin{displaymath}
\begin{array}{|c|c|c|c|c|c|}\hline&&&&&\\[-2ex]
\mbox{type} & \lambda & X_{\lambda} & \bfs & L(\nu) & \mbox{codim}
\\[1ex]\hline&&&&&\\[-2ex]
5a & 0 & \basis{e} & 1 & \nu & 1 \\[1.5ex]
5b & 0^2 & \basis{e,Ne} & T &
\left(\begin{array}{cc}\nu_1&\nu_2\\1&\nu_1\end{array}\right) & 2\\[3ex]
5c & 00 & \basis{a,b} & R &
\left(\begin{array}{rr}\nu_1&-\nu_3\\\nu_3&\nu_2\end{array}\right) & 3\\[3ex]
5d & 00 & \basis{a_1,a_2} & I &
\left(\begin{array}{cc}\nu_1&\nu_3\\\nu_3&\nu_2\end{array}\right) & 3
\\[3ex]\hline&&&&&\\[-2ex]
6a & 0 & \basis{e} & 1 & 0 & 0\\[1.5ex]
6b & (0^+)^3 & \basis{e,Ne,N^2e} &
\left(\begin{array}{rrr}&&1\\&-1&\\1&&\end{array}\right) &
\left(\begin{array}{rrr}0&\nu&0\\1&0&\nu\\0&1&0\end{array}\right) & 1\\[4ex]
6c & 0^2 & \basis{e,f,Ne,Nf} &
\left(\begin{array}{rr}&J\\-J&\end{array}\right) &
\left(\begin{array}{rrrr}
\nu_2&\nu_3-\nu_4&\nu_1&0\\
\nu_3+\nu_4&-\nu_2&0&\nu_1\\
1&0&\nu_2&\nu_3-\nu_4\\
0&1&\nu_3+\nu_4&-\nu_2
\end{array}\right) & 4\\[5ex]
6d & (\pm i\beta)^2 & \basis{e,f,Ne,Nf} &
\left(\begin{array}{rr}&J\\-J&\end{array}\right) &
\left(\begin{array}{rrrr}
0&-\beta-\nu_1&\nu_2&0\\
\beta+\nu_1&0&0&\nu_2\\
1&0&0&-\beta-\nu_1\\
0&1&\beta+\nu_1&0
\end{array}\right) & 2\\[5ex]
6e & (\pm i\beta)(\pm i\beta) & \basis{e_1,f_1,e_2,f_2} &
\left(\begin{array}{cc}I_2&\\&I_2\end{array}\right) &
\left(\begin{array}{rrrr}
0&-\beta-\nu_1&\nu_3&-\nu_4\\
\beta+\nu_1&0&\nu_4&\nu_3\\
-\nu_3&-\nu_4&0&-\beta-\nu_2\\
\nu_4&-\nu_3&\beta+\nu_2&0
\end{array}\right) &4\\[5ex]\hline&&&&&\\[-2ex]
7a & 0 & \basis{e, f} & J &
\left(\begin{array}{cc}\nu&0\\0&\nu\end{array}\right) & 1\\[3ex]
7b & 0^2 & \basis{e, f, Ne, N f} &
\left(\begin{array}{cc}&J\\J&\end{array}\right) &
\left(\begin{array}{rrrr}
\nu_1&0&\nu_2&\nu_3-\nu_4\\
0&\nu_1&\nu_3+\nu_4&-\nu_2\\
1&0&\nu_1&0\\
0&1&0&\nu_1
\end{array}\right) & 4\\[5ex]\hline&&&&&\\[-2ex]
8a & 0 & \basis{e, f} & J &
\left(\begin{array}{rr}
\nu_1&\nu_2-\nu_3\\\nu_2+\nu_3&-\nu_1
\end{array}\right) & 3\\[3ex]
8b & 0^2 & \basis{e,Ne} & J &
\left(\begin{array}{cc}0&\nu\\1&0\end{array}\right) & 1\\[3ex]
8c & (\pm i\beta)^2 & \basis{e,f,Ne,Nf} &
\left(\begin{array}{rr}&I\\-I&\end{array}\right) &
\left(\begin{array}{rrrr}
0&-\beta-\nu_1&\nu_2&0\\
\beta+\nu_1&0&0&\nu_2\\
1&0&0&-\beta-\nu_1\\
0&1&\beta+\nu_1&0
\end{array}\right) & 2\\[5ex]
8d & (\pm i\beta)(\pm i\beta) & \basis{e_1,f_1,e_2,f_2} &
\left(\begin{array}{cc}J&\\&J\end{array}\right) &
\left(\begin{array}{rrrr}
0&-\beta-\nu_1&\nu_3&-\nu_4\\
\beta+\nu_1&0&\nu_4&\nu_3\\
-\nu_3&-\nu_4&0&-\beta-\nu_2\\
\nu_4&-\nu_3&\beta+\nu_2&0
\end{array}\right) &4\\[5ex]\hline
\end{array}
\end{displaymath}
\caption{\textit{Unfoldings in eigenspaces of anti-automorphisms. See the
caption of Table \ref{tab:unfo1234} for an explanation of the
notation.\label{tab:unfo5678}}}
\end{table}

There are several examples where we have passing or splitting depending on
the signs. In $\bfs$-equivariant maps of type $\fR$, type 1 in Table
\ref{tab:types}, passing of real eigenvalues with different signs
and splitting of real eigenvalues with equal signs are both
codimension one phenomena. This can be inferred from 1c and 1d
in Table \ref{tab:unfo1234}. In infinitesimally $\bfs$-reversible maps of
type $\fR$, type 2 in Table \ref{tab:types}, only zero eigenvalues are
signed. At collisions real and imaginary eigenvalues generically
split. Maps of type 3 generically do not have eigenvalues on the real or
imaginary axis. Maps of type 4 generically do have opposite pairs of real
eigenvalues, but since there are no signs they split at collisions. Real
eigenvalues of maps of type 5 with equal signs pass but split when the
signs are different at collisions. This follows from 5b, 5c and 5d in Table
\ref{tab:unfo5678}. Similarly imaginary eigenvalues of anti-symmetric maps
of type 6 split or pass when the signs are different or equal
respectively. See 6e and 6f in Table \ref{tab:unfo5678}. Note that passing
is a codimension 3 phenomenon, but splitting is a codimension 1 phenomenon.
This should come as no surprise because the anti-symmetric maps are closely
related to infinitesimally symplectic maps of type 8. Here we have the same
codimensions for splitting and passing, see 8c and 8d in Table
\ref{tab:unfo5678}. Maps of type 7 can generically have real
eigenvalues. Since there are no signs they generically split at collisions.

%
\section{Generalizations}\label{sec:gens}
Here we will generalize the results for a single (anti)-automorphism of
order two to an abelian group $\Gamma$ of (anti)-automorphisms of order
two. In general the subset $\g$ in the Introduction will be an isotypic
component of the action of $\Gamma$ on $\gl(V)$, but for abelian groups
these are equivalent to simultaneous eigenspaces. To make this more
precise, let the $V$ be a finite dimensional real vector space and let
$\Gamma$ be a abelian group of (anti)-automorphisms of order two acting on
$\gl(V)$. Suppose $\Gamma$ is generated by $\gen{\gamma_1,\ldots,\gamma_p}$
with $\gamma_i^2=I$ for $i=1,\ldots,p$. Then the simultaneous eigenspaces
are given by
\begin{displaymath}
\gl_{\mu_1,\ldots,\mu_p}(V) = \{L \in \gl(V) \;|\; \gamma_1(L) =
\mu_1 L, \ldots, \gamma_p(L) = \mu_p L\},
\end{displaymath}
where the eigenvalues $\mu_i$ are $\pm 1$. The structure map associated to
$\gamma_i$ is denoted by $\bfs_i$. Apart from Example \ref{ex:motiv3} we
encounter this situation with infinitesimally reversible
equivariant and infinitesimally symplectic reversible equivariant maps. See
\cite{hlr2} for applications of the results of this article.

The theory developed for a single (anti)-automorphism almost immediately
extends to an abelian group of (anti)-automorphisms. Let us review Section
\ref{sec:nfeig} and make some comments. The structure maps can again be
taken orthogonal, but here we need to take a closer look at
anti-automorphisms, see Appendix \ref{sec:ortho}. From the proof of
Proposition \ref{pro:autjcd} it follows immediately that the
Jordan-Chevalley decomposition also holds in
$\gl_{\mu_1,\ldots,\mu_p}(V)$. Moreover the \emph{structure preserving}
transformation group $\G$ is the intersection of the structure preserving
transformations groups for each (anti)-automorphism $\gamma_i$. Once we
have identified the transformation group we can classify its orbits in
$\gl_{\mu_1,\ldots,\mu_p}(V)$. There is a Reduction Lemma similar to
\ref{lem:autreduc} where indecomposable $\bfs$-invariant subspaces are
replaced by indecomposable $\bfs_1,\ldots,\bfs_p$-invariant
subspaces. In the same way we have an Unfolding Lemma and an Orbit Splitting
Theorem for $\Gamma$. But in the latter we now have at most $2^p$
inequivalent $\G$-orbits.
\begin{remark}{\rm
The indecomposable $\bfs_1,\ldots,\bfs_p$-invariant subspaces can be
relatively large. Let us look at an $\fH$-linear map on $\fR^4$. The
quaternionic structure on $\fR^4$ is determined by two structure maps
$\bfc$ and $\bfq$ with $\bfc^2=-I$, $\bfq^2=-I$ and $\bfc \bfq=-\bfq \bfc$,
see \cite{hlr2}. Then $\fH$-linear maps on $\fR^4$ are defined as
$\gl_{1,1}(\fR^4) \defdas \{A \in \gl(\fR^4) \;|\; \phi_{\bfc}(A)=A,\;
\phi_{\bfq}(A)=A\}$, where $\phi_{\bfc}(A) \defdas \bfc^{-1}A\bfc$ and
$\phi_{\bfq}(A) \defdas \bfq^{-1}A\bfq$. Here $\Gamma = \{I, \phi_{\bfc},
\phi_{\bfq}, \phi_{\bfc} \circ \phi_{\bfq}\}$ since clearly $\phi_{\bfc}$
and $\phi_{\bfq}$ commute. Let $L \in \gl_{-1,-1}(\fR^4)$ have a real
eigenvalue $\alpha$, then the indecomposable $L,\bfc,\bfq$-invariant space
is $X_{\alpha} = V_{\alpha} \oplus \bfc V_{\alpha} \oplus \bfq V_{\alpha}
\oplus \bfc \bfq V_{\alpha}$.\hfill $\rhd$}
\end{remark}
\begin{example}{\rm
The maximum number of $2^p$ inequivalent $\G$-orbits occurs in an example
of an infinitesimally reversible symplectic linear map. We have already
encountered such maps in Example \ref{ex:motiv3}. They are elements of the
simultaneous eigenspace $\gl_{-1,-1}(\fR^{2n}) \defdas \{A \in
\gl(\fR^{2n}) \;|\; \phi_R(A) = -A,\; \psi_J(A) = -A\}$, where $\phi_R(A)
\defdas R^{-1}AR$ and $\psi_J(A) \defdas J^{-1} A^* J$. Here $\Gamma=\{id,
\phi_R, \psi_J, \psi_{RJ}\}$ is generated by $\phi_R$ and $\psi_J$, where
$\psi_J \circ \phi_R = \psi_{RJ}$. Let $L$ be a map in $\gl_{-1,-1}(\fR^4)$
with two blocks of double zero eigenvalues and a nilpotent part of height
two, see \cite{hov}. Then there are $4=2^2$ inequivalent $\G$-orbits in
$\gl_{-1,-1}(\fR^4)$. \hfill $\rhd$}
\end{example}

As mentioned in the Introduction our main motivation for studying
eigenspaces of (anti)-automorphisms of order two comes from real ordinary
differential equations. Other obvious generalizations apart from the one
given in this section are not necessarily in this context. Such
generalizations include single (anti)-automorphisms of finite order. Then
the eigenspace $\gl_{\mu}(V)$ need not be real. This problem can be
overcome by taking the real invariant space $\gl_{\mu,\bar{\mu}}(V)$
as the object of study, though one could also look at the complex space
$\gl_{\mu}(V)$ in its own right.
One could also look at abelian and non-abelian groups
generated by (anti)-automorphisms of finite order, and more generally still
one might consider general compact groups of (anti)-automorphisms.
We will not pursue these matters here.

%
\section*{Acknowledgments}
It is a great pleasure to thank Richard Cushman for valuable comments and
discussions. This research was supported by the UK Engineering and Physical
Sciences Research Council EPSRC (IH, JSWL), the Nuffield Foundation
(JSWL), and by European Community funding for the Research Training Network
`MASIE' (HPRN-CT-2000-00113).

\appendix
%
\section{Standard Form of a Bilinear Form}\label{sec:gsbc}
Here we generalize a result of Burgoyne \& Cushman \cite{bc1,bc2}, which
in turn is
based on a theorem of Springer \& Steinberg \cite{ss}, to obtain a normal
form for a skew or indefinite symmetric bilinear form on $V$ which respects
the splitting of $V$ in the Reduction Lemma \ref{lem:autreduc}.

Let $V$ be real vector space and let $\omega$ be a nondegenerate bilinear
form on $V$ which is either symmetric or skew, that is for each nonzero $x
\in V$ there exists an $y \in V$ such that $\omega(x,y) \neq 0$ and for all
$x,y \in V$, $\omega(y,x) = \eps \omega(x,y)$ with $\eps = \pm 1$. If
$\inprod{\cdot}{\cdot}$ is the standard inner product on $V$ then there is
an invertible linear map $\bfs$ with $\bfs^* = \eps \bfs$ such that
$\omega(x,y) = \inprod{x}{\bfs y}$. We may assume that after scaling $\bfs^2
= \eps I$, so $\bfs$ is orthogonal.

Let $L$ be a linear map on $V$ such that for all $x,y \in V$, $\omega(Lx,y)
= \mu \omega(x,Ly)$ with $\mu = \pm 1$. This is equivalent to
$\inprod{Lx}{\bfs y} = \mu \inprod{x}{\bfs Ly}$ or $L^*\bfs = \mu \bfs
L$. Let $\psi_{\bfs}(L) = \bfs^{-1}L^*\bfs$. Then $L$ satisfies
$\psi_{\bfs}(L) = \mu L$.

Now we assume that $V$ is an indecomposable $L$,$\bfs$-invariant space. If
$L=S+N$ is the Jordan-Chevalley decomposition of $L$ then there is an
$S$-invariant complement $W$ of $NV$ in $V$ such that $V=W \oplus NW \oplus
\cdots \oplus N^{n-1}W$, where $n$ is the height of $N$. Although $\bfs$
has a normal form as in Lemma \ref{pro:sprops}, since $\bfs$ is
orthogonal, we wish to find a normal form of $\bfs$ which respects the
above splitting of $V$.

If $\omega$ is definite then it must be symmetric and thus $\bfs = \pm
I$. Then every $L$ satisfying $\psi_{\bfs}(L) = \mu L$ is
semi-simple. Therefore in the present situation we need only consider
indefinite forms. The main result of this section can now be stated.
\begin{proposition}\label{pro:gsbc}
Let $\omega$ and $L$ be defined as above. For every $S$-invariant
complement $W_1$ of $NV$ in $V$ there is an invertible transformation $g$
such that $W=gW_1$ is again $S$-invariant and on $W \oplus NW \oplus \cdots
\oplus N^{n-1}W$ the matrix of $\bfs$ takes the form
\begin{displaymath}
\left(\begin{array}{ccc}&&\ast\\&\adots&\\\ast&&\end{array}\right),
\end{displaymath}
where $\ast$ is a $m \times m$ block and $m = \dim W$.
\end{proposition}

The procedure to transform $W_1$ runs as follows. Note that the matrix of
$\bfs$ is upper triangular with respect to the main
anti-diagonal. Furthermore blocks of $\bfs$ on an anti-diagonal differ only
by a sign, see the proof of Proposition \ref{pro:gsbc} below. By adding a
component of $NW_1$ to $W_1$ and setting $W_2 = W_1 + Ng_1W_1$ for a map
$g_1$
to be specified later, we clear the first co-anti-diagonal. Then we set $W_3 =
W_2 + N^2g_2W_2$, clearing the second co-anti-diagonal without affecting
the first. This process stops after $n-1$ steps.

The following lemmas are useful in the proof of Proposition
\ref{pro:gsbc}. Although there is freedom in choosing an $S$-invariant
complement $W$ of $NV$ in $V$ the space $N^{n-1}W$ is unique.
\begin{lemma}
If $W$ is an $S$-invariant complement $W$ of $NV$ in $V$, then $N^{n-1}W$
is unique.
\end{lemma}

\textbf{Proof.} Let $n$ be the height of $N$ on $V$. For every $x \in
N^{n-1}W$ we have $Nx=0$. Since the eigenspaces of $N$ are unique
$N^{n-1}W$ is unique. Thus if $U$ is also an $S$-invariant complement $W$
of $NV$ in $V$ then $N^{n-1}U = N^{n-1}W$.\hfill $\Box$

The matrix of $\bfs$ with respect to a basis in $V=W \oplus NW \oplus \cdots
\oplus N^{n-1}W$ has the following properties.
\begin{lemma}
Let $\omega(x,y)=\inprod{x}{\bfs y}_V$ and blocks of $\bfs$ are denoted by
$\beta_{ij}$. Then $\beta_{i,j} = \mu \beta_{i-1,j+1}$ and $\beta_{ij}=0$ if
$i+j \geq n$.
\end{lemma}
\textbf{Proof.} This follows immediately from $\omega(N^ix, N^jy) =
\mu\omega(N^{i-1}x, N^{j+1}y) = \mu^i \omega(x,N^{i+j}y)$.\hfill $\Box$

In the proof of the proposition we need the bilinear forms $\tau_j(x,y) =
\omega(N^jx,y)$ on $W$ for $j=1,\ldots,n-1$. Furthermore let $T_j$ be
defined as $\tau_j(x,y) = \inprod{x}{T_jy}_W$ for all $x,y \in W$.
\begin{lemma}
$T_{n-1}$ is an invertible map.
\end{lemma}
\textbf{Proof.} This follows from the fact that $\tau_{n-1}$ is
nondegenerate on $W$. The form $\omega$ is nondegenerate on $V$, so for every
nonzero $x \in V$ there is a $y \in V$ such that $\omega(x,y) \neq 0$. In
particular for every $x \in W$ there is a $y \in V$ such that
$\omega(N^{n-1}x,y) \neq 0$. Every such $y$ has a unique decomposition
$y=y_1+y_2$ with $y_1 \in W$ and $y_2 \in NV$. Then $0 \neq
\omega(N^{n-1}x,y) = \mu^{n-1}\omega(x,N^{n-1}y) =
\mu^{n-1}\omega(x,N^{n-1}y_1 + N^{n-1}y_2) =
\mu^{n-1}\omega(x,N^{n-1}y_1)$. Thus $\tau_{n-1}$ is nondegenerate on $W$
and therefore $T_{n-1}$ is invertible on $W$.\hfill $\Box$

\textbf{Proof of Proposition \ref{pro:gsbc}.} Note that the blocks
$\beta_{ij}$ differ from the matrices of $T_{j-1}$ by a sign
only. Therefore if the bilinear forms $\tau_j$ are identically zero on $W$
for $j=1,n-2$, the matrix of $\bfs$ has the desired form.

Suppose $\tau_{n-2} \neq 0$ on $W$. Let $W_2=W_1+NgW_1$ where $g$ is chosen
so that $\tau_{n-2}((I+Ng)x, (I+Ng)y) = 0$ for all $x,y \in W_1$. After a
short computation we find $\tau_{n-2}(x+Ngx, y+Ngy) = \inprod{x}{(T_{n-2} +
\mu T_{n-1}g + g^*T_{n-1})y}_W$. Assuming that $g^*T_{n-1}=\mu T_{n-1}g$ we
set $g=-\halfje \mu T_{n-1}^{-1} T_{n-2}$. Then $\tau_{n-2} = 0$ on $W_2$.

Now assume that $\tau_{n-2}=0,\ldots,\tau_{n-j+1}=0$ on $W_{j-1}$. Let $W_j
= W_{j-1} + N^{j-1}gW_{j-1}$, where $g$ is chosen so that
$\tau_{n-j}((I+N^{j-1}g)x, (I+N^{j-1}g)y) = 0$ for all $x,y \in
W_{j-1}$. Again we find $\tau_{n-j}(x+N^{j-1}gx, y+N^{j-1}gy) =
\inprod{x}{(T_{n-j} + \mu^{j-1} T_{n-1}g + g^*T_{n-1})y}_W$ and we set
$g=-\halfje \mu^{j-1} T_{n-1}^{-1} T_{n-j}$, so that $\tau_{n-j} = 0$ on
$W_j$. It is easily checked that now $\tau_{n-2}=0,\ldots,\tau_{n-j}=0$ on
$W_j$.

We still have to check that $g^*T_{n-1} = \mu^{j-1} T_{n-1} g$ in each
step, but this follows from $T_{n-j}^*T_{n-1} = \mu^{j-1} T_{n-1} T_{n-j}$
because $N^*T_{n-1}=\mu T_{n-1}N$. Furthermore it is easy to see that each
$W_j$ is $S$-invariant.\hfill $\Box$

%
\section{Orthogonality of Structure Maps}\label{sec:ortho}
In this appendix we give a precise statement of the properties of structure
maps associated to the generators of an abelian group $\Gamma$ of
(anti)-automorphisms of order two. Their properties are essentially those
of a structure map of a single (anti)-automorphism, but it is not a priori
clear that we can transform them as in Proposition \ref{pro:sprops} so
that they can all be assumed to be orthogonal. It might happen that a
transformation which takes one structure map in good shape spoils
another. The proposition below shows that this does not happen because of
the commutation relations of the (anti)-automorphisms.
\begin{proposition}\label{pro:gsprops}
Let $\Gamma$ be a abelian group of (anti)-automorphisms of order two on
$\gl(V)$ generated by $\gen{\gamma_1,\ldots,\gamma_p}$. Then we may assume
that $\gamma_i = \gamma_{\bfs_i}$ where the structure maps $\bfs_i$ have
the following properties:
\begin{enumerate}\abc\parskip 0pt \itemsep 0pt
\item $\det \bfs_i = \pm 1$.
\item $\bfs_i$ is orthogonal.
\item $\bfs_i^2 = \pm I$.
\item $\bfs_i \bfs_j = \pm \bfs_j \bfs_i$.
\end{enumerate}
\end{proposition}

\textbf{Proof.} Throughout the proof $\bfs$ and $\bft$ will be any
pair of structure maps from the set $\{\bfs_1 \ldots \bfs_p\}$.

Part a) follows from $\gamma_{\alpha \bfs} =\gamma_{\bfs}$
for all $\alpha \in \fR$, so we can scale $\bfs$ such that $\det \bfs = \pm
1$. For every pair $\gamma_{\bfs}, \gamma_{\bft} \in \Gamma$ we have
$\gamma_{\bfs} \circ \gamma_{\bft} = \gamma_{\bft} \circ \gamma_{\bfs}$ and
$\gamma_{\bfs}^2=I$, $\gamma_{\bft}^2=I$. We distinguish three different
cases.
\begin{enumerate}
\item $\gamma_{\bfs} = \phi_{\bfs}$ and $\gamma_{\bft} = \phi_{\bft}$ are
automorphisms. From $\phi_{\bfs}^2 = I$ we have $\bfs^2 = \pm I$, so $\bfs$
is semi-simple and has eigenvalues $\pm 1$ or $\pm i$. In order that
$\phi_{\bfs} \circ \phi_{\bft} = \phi_{\bft} \circ \phi_{\bfs}$ we must
have $\bfs \bft = \pm \bft \bfs$. Then $\gen{\bfs, \bft}$ generates a
finite group. By a transformation, corresponding to averaging the inner
product on $V$ over this group, we obtain that $\bfs$ and $\bft$ are
orthogonal. We can do this at once for all structure maps associated to
automorphisms in $\gen{\gamma_1,\ldots,\gamma_p}$.
\item $\gamma_{\bfs} = \phi_{\bfs}$ is an automorphism and $\gamma_{\bft} =
\psi_{\bft}$ is an anti-automorphism. Because of 1) we assume that $\bfs$
has properties a), b) and c). Now $\psi_{\bft}^2=I$ implies $\bft^*=\pm
\bft$ and so $\bft$ is semi-simple, moreover $\bft$ has either real or
purely imaginary eigenvalues. Then it follows from $\phi_{\bfs} \circ
\psi_{\bft} = \psi_{\bft} \circ \phi_{\bfs}$ that $\bfs \bft = \pm \bft
\bfs$. The latter implies that there exist indecomposable $\bfs,
\bft$-invariant subspaces on which $\bft$ has either real eigenvalues in
configurations $\alpha$, $\alpha \alpha$ or $\pm\alpha$, or purely
imaginary eigenvalues in configurations $\pm i\beta$ or $(\pm i\beta) (\pm
i\beta)$. A scaling transformation acts on this subspace as $g=\rho I$,
taking the eigenvalues of $\bft$ to $\pm 1$ or $\pm i$. Since $\bfs$
transforms as $\bfs \mapsto \bfs g^{-1}$ it is invariant under the scaling
$g$. Thus $\bft$ has properties a), b), c) and d).
\item $\gamma_{\bfs} = \psi_{\bfs}$ and $\gamma_{\bft} = \psi_{\bft}$ are
anti-automorphisms. Because of 2) we assume that $\bfs$ has properties a),
b) and c). Again $\psi_{\bft}^2=I$ implies $\bft^*=\pm \bft$ and so $\bft$
is semi-simple, moreover $\bft$ has either real or purely imaginary
eigenvalues. From $\psi_{\bfs} \circ \psi_{\bft} = \psi_{\bft} \circ
\psi_{\bfs}$ we infer that $\bft \bfs = \pm \bfs \bft^{-1}$. Again we look
for a transformation that takes eigenvalues of $\bft$ to $\pm 1$ or $\pm i$
but leaves $\bfs$ invariant. Summarising we have eight different cases
$\bfs^2 = \eps_1 I$, $\bfs^* = \eps_1 \bfs$, $\bft^* = \eps_2 \bfs$, $\bft
\bfs = \eps_3 \bfs \bft^{-1}$, where $\eps_i = \pm 1$. Let us look at
$\eps_1=\eps_2=-1$ and $\eps_3=\eps=\pm 1$, the other cases being very
similar. On any indecomposable $\bfs, \bft$-invariant subspace $\bft$ has
eigenvalues $\pm i\beta$ and $\pm \frac{i}{\beta}$, $\beta > 0$. Suppose
$e$ is a vector such that $\bft^2 e=-\beta^2 e$. Let $f$ be defined as $f
\defdas \frac{1}{\beta}\bft e$. Then $\bft \bfs e = \eps \bfs \bft^{-1} e =
-\eps \frac{1}{\beta} \bfs f$ and $\bft \bfs f = \eps \bfs \bft^{-1} f =
\eps \frac{1}{\beta} \bfs e$. Since $\bfs^2=-I$, $\bfs^2 v = -v$ for each
vector $v$. Thus on the basis $\basis{e,f,\bfs e,\bfs f}$, $\bft$ and
$\bfs$ have the following matrices
\begin{displaymath}
\bft=\left(\begin{array}{rrrr}0&-\beta&&\\\beta&0&&\\
&&0&\eps\beta^{-1}\\&&-\eps\beta^{-1}&0\end{array}\right),\quad
\bfs=\left(\begin{array}{rrrr}&&-1&0\\&&0&-1\\
1&0&&\\0&1&&\end{array}\right).
\end{displaymath}
Now let $g$ be a transformation with blockdiagonal matrix
$\diag{a^{-1}I_2,aI_2}$ with $a=\sqrt{\beta}$. Then $\bft$ and $\bfs$
transform as
\begin{displaymath}
g\bft g^* = \left(\begin{array}{rrrr}0&-1&&\\1&0&&\\
&&0&\eps\\&&-\eps&0\end{array}\right),\quad
g\bfs g^* = \left(\begin{array}{rrrr}&&-1&0\\&&0&-1\\
1&0&&\\0&1&&\end{array}\right) = \bfs.
\end{displaymath}
So we see that $\bft$ has the properties listed in the lemma.\hfill $\Box$
\end{enumerate}


%


\begin{thebibliography}{99}

\bibitem{arn1} V.I. Arnol'd, ``Geometrical Methods in the Theory of Ordinary
Differential Equations'', Springer-Verlag, New York, 1983.

\bibitem{bre} G.E. Bredon, ``Introduction to Compact Transformation
Groups'', Pure and Applied Mathematics, vol.46, Academic Press, New York,
1972.

\bibitem{bdst} H.W. Broer, F. Dumortier, S. van Strien and F. Takens,
``Structures in Dynamics'', Studies in Mathematical Physics \textbf{2}
(E. van Groesen and E.M. de Jager eds.), North-Holland, Amsterdam, 1991.

\bibitem{bc1} N. Burgoyne and R.H. Cushman, The decomposition of a linear
mapping, \textit{Linear Algebra and its Applications} \textbf{8} (1974),
515-519.

\bibitem{bc2} N. Burgoyne and R.H. Cushman, Normal forms for real linear
Hamiltonian systems, in: ``The 1976 Ames Research Center (NASA) Conference on
Geometric Control Theory'', (C. Martin and R. Hermann, Eds), pp. 483-529,
Math. Sci. Press., Brookline, Mass., 1977.

\bibitem{bc3} N. Burgoyne and R.H. Cushman, Conjugacy classes in linear
groups, \textit{J. Alg.} {\bf 44} (1977), 333-362.

\bibitem{cot} C. Cotter, ``The 1:1 Semi-Simple Resonance'', Ph.D. thesis,
University of California at Santa Cruz, 1986.

\bibitem{cs} R.H. Cushman and J.A. Sanders, Nilpotent normal forms and
representation theory of $\sl(2,\fR)$, in: ``Multiparameter Bifurcation
Theory'', eds. M. Golubitsky and J.M. Guckenheimer, Contemporary Mathematics
Vol. \textbf{56}, AMS 1985, 31-51.

\bibitem{dem} U. Dempwolff, Normal forms and fixed subspaces of semi-linear
maps, \textit{Bolletino U.M.I.} \textbf{7} 4-A (1990), 209-218.

\bibitem{dpwz} D.Z. Djukovic, J. Patera, P. Winternitz and H. Zassenhaus,
Normal forms of elements of classical real and complex Lie and
Jordan-algebras, \textit{J. Math. Phys.} \textbf{24} No. 6 (1983), 1363-1374.

\bibitem{fv} H. Freudenthal and H. de Vries, ``Linear Lie Groups'', Academic
Press, New York, 1969.

\bibitem{gal} D.M. Galin, Versal deformations of linear Hamiltonian
systems, \textit{Amer. Math. Soc. Transl.} \textbf{118} No. 2 (1982), 1-12.

\bibitem{gan} F.R. Gantmacher, ``Theory of Matrices'', Chelsea, New York,
1959.

\bibitem{hov} I. Hoveijn, Versal Deformations and normal forms for
reversible and Hamiltonian linear systems, \textit{J. Diff. Eq.}
\textbf{126} No. 2 (1996), 408-442.

\bibitem{hlr2} I. Hoveijn, J.S.W. Lamb and R.M. Roberts, Reversible
equivariant linear systems: normal forms and unfoldings, in preparation.



\bibitem{hum} J.E. Humphreys, ``Introduction to Lie Algebras and
Representation Theory'', GTM 9, Springer-Verlag, New York, Berlin, 1972.

\bibitem{ioo} G. Iooss, A codimension 2 bifurcation for reversible vector
fields, in: ``Proceedings Normal forms and Homoclinic
Chaos''. eds. W. Langford and W. Nagata, Fields Institute Communications
\textbf{4} (1995), 201-218.

\bibitem{jac} N. Jacobson, Pseudo-linear transformations, \textit{Annals of
Mathematics} \textbf{38} No. 2 (1937), 484-507.

\bibitem{koc} H. Ko\c{c}ak, Normal forms and versal deformations of linear
Hamiltonian systems. \textit{J. Differential Equations} \textbf{51} (1984),
359-407.

\bibitem{lr} J.S.W. Lamb and R.M. Roberts, Reversible equivariant linear
systems, \textit{Journal of Differential Equations} \textbf{159} (1999),
239-279.

\bibitem{mee} J.C. van der Meer, ``The Hamiltonian Hopf Bifurcation'',
Lecture Notes in Mathematics 1160, Springer-Verlag, New York, Berlin, 1985.

\bibitem{mel} I. Melbourne, Versal unfoldings of equivariant linear
Hamiltonian vector fields, \textit{Math. Proc. Cambridge Philos. Soc.}
\textbf{114} (1993), 559-573.

\bibitem{md} I. Melbourne and M. Dellnitz, Normal forms for linear
Hamiltonian vector fields commuting with the action of a compact Lie group,
\textit{Math. Proc. Cambridge. Philos. Soc.} \textbf{114} (1993), 235-268.

\bibitem{pal} K.J. Palmer, Linearisation of reversible systems,
\textit{J. Math. Anal. Appl.} \textbf{60} (1977), 794.

\bibitem{pr2} J. Patera and C. Rousseau, Versal defomations of elements of
classical Jordan algebras, \textit{J. Math. Phys.} \textbf{24} No. 6
(1983), 1363--1374.

\bibitem{sev1} M.B. Sevryuk, Reversible linear systems and their versal
deformations, \textit{J. Sov. Math.} \textbf{60} No. 5 (1992), 1663-1680.

\bibitem{shi} C.W. Shih, Normal forms and versal deformations of linear
involutive dynamical systems, \textit{Chinese J. Math.} \textbf{21}, No. 4
(1993), 333-347.

\bibitem{ss} T.A. Springer and R. Steinberg, ``Conjugacy Classes'', Lecture
Notes in Mathematics, Vol. 131, Springer, Berlin, 1970.

\bibitem{var} V.S. Varadarajan, ``Lie Groups, Lie Algebras, and their
Representations'', GTM 102, Springer, 1984.

\bibitem{wil} J. Williamson, On the algebraic problem concerning the normal
forms of linear dynamical systems, \textit{Amer. J. Math.} \textbf{58} (1936),
141-163.

\bibitem{wie} N.A. Wiegman, Some theorems on matrices with real
quaternionic elements, \textit{Can. J. of Math.} (1955).

\end{thebibliography}
\end{document}